\documentclass{article}
\usepackage[english]{babel}

\usepackage{latexsym}
\usepackage{amsmath}
\usepackage{amsfonts}
\usepackage{mathrsfs}
\usepackage{amstext}
\usepackage{amssymb}
\usepackage{amsthm}       

\usepackage{xspace}

\newtheorem{thm}{Theorem}[section]
\newtheorem{prop}[thm]{Proposition}
\newtheorem{cor}[thm]{Corollary}

\newtheorem{lem}[thm]{Lemma}

\begin{document}

\title{A real analytic family of fundamental solutions of elliptic partial differential operators with real constant coefficients}
\author{M.~Dalla Riva}
\maketitle

\section{Introduction} 
We fix once for all
\[
n,k\in\mathbb{N}\,,\quad n\ge 2\,,\quad k\ge 1\,.
\] Here $\mathbb{N}$ denotes the set of natural numbers including $0$. Then we denote by $N(2k,n)$ the set of all multi-indexes $\alpha\equiv(\alpha_1,\dots,\alpha_n)\in\mathbb{N}^n$ such that $|\alpha|\equiv\alpha_1+\dots+\alpha_n\leq 2k$. We denote by  $R(2k,n)$ the set of the functions $\mathbf{a}\equiv(a_{\alpha})_{\alpha\in N(2k,n)}$  from $N(2k,n)$ to $\mathbb{R}$. We note that $R(2k,n)$ can be identified with a finite dimensional real vector space. Accordingly we endow $R(2k,n)$ with the corresponding Euclidean norm $|\mathbf{a}|\equiv(\sum_{\alpha\in N(2k,n)} a_\alpha^2)^{1/2}$. Then we set 
\[
\begin{split}
&\mathscr{E}_R(2k,n)\\
&\quad\equiv\biggl\{\mathbf{a}=(a_{\alpha})_{\alpha\in N(2k,n)}\in R(2k,n)\,:\,\sum_{|\alpha|=2k}a_\alpha\xi^\alpha\neq 0\text{ for all }\xi\in\mathbb{R}^n\setminus\{0\}\biggr\}\,.
\end{split}
\]  We note that $\mathscr{E}_R(2k,n)$ is an open non-empty subset of $R(2k,n)$.   Then, for each  $\mathbf{a}\in \mathscr{E}_R(2k,n)$ we denote by $\mathbf{L}[\mathbf{a}]$ the partial differential operator  defined by
\[
\mathbf{L}[\mathbf{a}]\equiv \sum_{\alpha\in N(2k,n)}a_{\alpha}\partial_{x_1}^{\alpha_1}\dots\partial_{x_n}^{\alpha_n}
\] (see also Section \ref{Notation} below). So that,  $\mathbf{L}[\mathbf{a}]$ is a real constant coefficients elliptic partial differential operator on $\mathbb{R}^n$ of order  $2k$.  

The aim of this paper is to show the construction  of a function $S$  such that
\begin{eqnarray}\label{Sanal}
&& S\text{ is a real analytic function from }\mathscr{E}_R(2k,n)\times(\mathbb{R}^n\setminus\{0\})\text{ to }\mathbb{R}\,;\\
\label{Sfs}
&& S(\mathbf{a},\cdot)\text{ is a fundamental solution of }\mathbf{L}[\mathbf{a}]\text{ for all }\mathbf{a}\in\mathscr{E}_R(2k,n)\,.
\end{eqnarray}
Condition \eqref{Sfs} means that $S(\mathbf{a},\cdot)$ defines a locally integrable function on $\mathbb{R}^n$ such that $\mathbf{L}[\mathbf{a}]S(\mathbf{a},\cdot)=\delta_0$ in the sense of distributions on $\mathbb{R}^n$, where $\delta_0$ denotes the delta Dirac distribution with mass at $0$.  Also, we wish to investigate regularity properties in the frame of Schauder spaces and jump properties of the single layer potential corresponding to the fundamental solution $S(\mathbf{a},\cdot)$. 

In Theorem \ref{Sa} below we introduce a function $S$ which satisfies the conditions in \eqref{Sanal} and \eqref{Sfs}.  Then, in Theorem \ref{fsthm1} we provide a suitably detailed expression for  $S$. In particular, we show that there exist a real analytic function $A$ from $\mathscr{E}_R(2k,n)\times\partial\mathbb{B}_{n}\times\mathbb{R}$ to $\mathbb{R}$,  and real analytic functions $B$ and $C$ from $\mathscr{E}_R(2k,n)\times\mathbb{R}^n$ to $\mathbb{R}$ such that 
\begin{equation}\label{Sform}
S(\mathbf{a},x)= |x|^{2k-n}A(\mathbf{a},x/|x|,|x|)+ \log |x|\,  B(\mathbf{a},x) +C(\mathbf{a},x)
\end{equation} for all $(\mathbf{a},x)\in\mathscr{E}_R(2k,n)\times(\mathbb{R}^n\setminus\{0\})$, where $B$ and $C$ are identically $0$ if the dimension $n$ is odd. The functions $A$ and $B$ play an important role when we consider the regularity and jump properties of the single layer potential corresponding to $S(\mathbf{a},\cdot)$. Therefore, we investigate the power series expansion of $A(\mathbf{a},\theta,r) $ with respect to the ``radius'' variable $r$ and of $B(\mathbf{a},x)$ with respect to the spatial variable $x$ (cf.~Theorem \ref{fsthm1} below). Then we pass to consider the principal term $\mathbf{L}_0[\mathbf{a}]$ of $\mathbf{L}[\mathbf{a}]$. Namely, we set
\begin{equation}\label{L0}
\mathbf{L}_0[\mathbf{a}]\equiv\sum_{|\alpha|=2k}a_\alpha\partial_{x_1}^{\alpha_1}\dots\partial_{x_n}^{\alpha_n}
\end{equation} for all $\mathbf{a}\in\mathscr{E}_R(2k,n)$ (see also Section \ref{Notation} below). In Theorem \ref{fsthm2} we show that there exists a real analytic function $S_0$ from $\mathscr{E}_R(2k,n)\times(\mathbb{R}^n\setminus\{0\})$  to $\mathbb{R}$  such that 
$S_0(\mathbf{a},\cdot)$  is a fundamental solution of $\mathbf{L}_0[\mathbf{a}]$  for all $\mathbf{a}\in\mathscr{E}_R(2k,n)$. We also provide an expression for $S_0$ in terms of the coefficients of the power series expansion of the functions $A$ and  $B$ which appears in \eqref{Sform}.

Then we turn to consider the single layer potential corresponding to the fundamental solution $S(\mathbf{a},\cdot)$. To do so, we fix once for all
\[
m\in\mathbb{N}\setminus\{0\}\quad\text{and}\quad\lambda\in]0,1[\,.
\] Then we fix a set $\Omega$ such that
\[
\Omega\text{ is an open bounded subset of }\mathbb{R}^n\text{ of class }C^{m,\lambda}\,.
\]  For the definition of functions and sets of the usual Schauder class $C^{h,\lambda}$, with $h\in\mathbb{N}$, we refer for example to Gilbarg and Trudinger~\cite[\S 4.1 and \S 6.2]{GT01}.
 Then we denote by $v[\mathbf{a},\mu]$ the single layer potential with density $\mu\in C^{m-1,\lambda}(\partial\Omega)$ corresponding to the fundamental solution $S(\mathbf{a},\cdot)$, $\mathbf{a}\in\mathscr{E}_R(2k,n)$. Namely, $v[\mathbf{a},\mu]$ is the function from $\mathbb{R}^n$ to $\mathbb{R}$ defined by  
\begin{equation}\label{sl}
v[\mathbf{a},\mu](x)\equiv\int_{\partial\Omega} S(\mathbf{a},x-y)\mu(y)\,d\sigma_y\qquad\forall x\in\mathbb{R}^n\,,
\end{equation} for all $\mathbf{a}\in\mathscr{E}_R(2k,n)$ and all $\mu\in C^{m-1,\lambda}(\partial\Omega)$, where $\partial\Omega$ denotes the boundary of $\Omega$ and $d\sigma$ denotes the area element on $\partial\Omega$. In Theorem \ref{2k-2} we show that $v[\mathbf{a},\mu]$ is a function of class $C^{2k-2}$ on $\mathbb{R}^n$. In Theorem \ref{fs-vs-m} we show that the restriction of $v[\mathbf{a},\mu]$ to the closure of $\Omega$ is a function of class $C^{m+2k-2,\lambda}$ and the restriction of $v[\mathbf{a},\mu]$ to $\mathbb{R}^n\setminus\Omega$ belongs to $C^{m+2k-2,\lambda}$ in a local sense which will be clarified. However, the derivatives  of order $2k-1$ of $v[\mathbf{a},\mu](x)$ with respect to the spatial variable $x$  are not continuous on $\mathbb{R}^n$ and display a jump across the boundary  of $\Omega$. In Theorem \ref{jump} we describe such a jump property.

We observe that the construction of the function $S$ presented here is based on the construction of a fundamental solution provided by John in \cite[Chapter III]{J55}. Moreover, the result which we show in Theorem \ref{Sa} resemble those which were proved by Tr{\`e}ves in \cite{MR0149084} and by Mantlik in \cite{MR1136596, MR1107027}. We also note that Tr{\`e}ves and Mantlik consider more general assumptions on the operators. However,  our result is not a corollary and differs from those of . Indeed, the joint real analyticity of  $S(\mathbf{a},x)$ upon the variables $\mathbf{a}\in\mathscr{E}_R(2k,n)$ and $x\in\mathbb{R}^n\setminus\{0\}$ does not follow by Tr{\`e}ves and Mantlik results. 
Also, the suitably detailed expression which we provide for $S$ in Theorem \ref{fsthm1} cannot be deduced by Tr{\`e}ves and Mantlik results (see also equality \eqref{Sform}).

It is also worth noting that  in case of operators of order $2$, the real analytic function $S$ which we introduce here is of the type considered by Lanza de Cristoforis and the author in \cite{DaLa10}. Accordingly, one could verify real analyticity results for the (joint) dependence of the single layer potential $v[\mathbf{a},\mu]$ upon perturbations of the variables $\mathbf{a}$ and $\mu$ and upon perturbations of the ``shape'' of the support of integration $\partial\Omega$ (cf.~\cite[Theorem 5.6]{DaLa10}, see also \cite[Theorem 5.1]{Dal09}). Then, one could exploit such results to analyze the dependence of solutions of boundary value problems upon perturbations of the domain of definition and of the coefficients of the corresponding operators. This program has been carried out for boundary value problems for the Laplace operator by Lanza de Cristoforis (cf., {\it e.g.}, \cite{Lan05a, Lan05b}, see also Lanza de Cristoforis and Rossi \cite{LaRo04, LaRo08})  and for the equations of linearized elasticity (cf.~\cite{Dal07}). This paper can be considered as a first step to generalize such an approach to the case of general elliptic partial differential operators with real constant coefficients. Also, the construction of the function $S$ presented here can be extended to the case of partial differential operators with complex or quaternion constant coefficients (cf.~\cite{DaMu12}) and to the case of particular systems of differential operators (cf.~\cite{Dal09}).

The paper is organized as follows.  In Section \ref{Notation} we introduce some standard notation and we recall a classical result on real analytic functions. In Section \ref{sec2.1} we define the G\"unter tangential derivative  $\mathscr{D}_\theta$ on the boundary of the unit ball $\mathbb{B}_n$. Then we prove some standard properties of $\mathscr{D}_\theta$.  Sections \ref{vsec} and \ref{Wsec} are devoted to analysis of some auxiliary functions and are therefore rather technical. In  Section \ref{vsec}  we consider a function $v$ which is solution of the equation $\mathbf{L}[\mathbf{a}]v=1$ and which vanish together with its derivatives of order $\le 2k-1$ on a hyperplane of $\mathbb{R}^n$.  In Section \ref{Wsec}, we introduce  the auxiliary functions $W_0$, $W_1$, and $W_2$. We show that one can exploit $W_0$, $W_1$, and $W_2$ to define a distribution $S_\mathbf{a}$ which is a fundamental solution of the operator $\mathbf{L}[\mathbf{a}]$ (cf.~Proposition \ref{S}).   In Section \ref{Ssec} we are ready to introduce our functions $S$ and $S_0$ and  we prove our main Theorems \ref{Sa}, \ref{fsthm1}, and \ref{fsthm2}. We observe that in Theorem \ref{Sa} we verify that the distribution $S_\mathbf{a}$ of Section \ref{Wsec} coincides with the distribution  defined by the function $S(\mathbf{a},\cdot)$.  In the last Section \ref{slsec}, we consider the single layer potential  $v[\mathbf{a},\mu]$ and we prove Theorems \ref{2k-2}, \ref{fs-vs-m}, and \ref{jump} where we investigate regularity and jump properties of $v[\mathbf{a},\mu]$.

\par\medskip

\section{Some notation and preliminaries}\label{Notation}
For standard definitions of calculus in normed spaces, we refer, {\it e.g.}, to Cartan \cite{Car71} and to  Prodi and Ambrosetti \cite{PA73}. We understand that a finite product of normed spaces is equipped with the sup-norm of the norm of the components, while we use the euclidean norm for $\mathbb{R}^n$.

 For all $x\in\mathbb{R}^n$, $x_j$ denotes the $j$-th coordinate of $x$, $|x|$ denotes the euclidean modulus of $x$, and $\mathbb{B}_n$ denotes the unit ball $\{x\in\mathbb{R}^n\,:\,|x|<1\}$. A dot `$\cdot$' denotes the inner product in $\mathbb{R}^n$. 
If $\mathcal{X}$ is a subset of $\mathbb{R}^n$, then $\mathrm{cl}\mathcal{X}$ denotes the closure of $\mathcal{X}$ and $\partial \mathcal{X}$ denotes the boundary of $\mathcal{X}$. If $\mathcal{O}$ is an open subset of $\mathbb{R}^n$, and $f$ is a function from $\mathcal{O}$ to $\mathbb{R}$, and $x\in \mathcal{O}$, then the partial derivative of $f$ with respect to $x_j$ at $x$ is denoted by $\partial_{x_j} f(x)$. Then $\partial_x^\alpha f(x) \equiv \partial_{x_1}^{\alpha_1}\dots\partial_{x_n}^{\alpha_n} f(x)$ for all multi-index $\alpha\equiv(\alpha_1,\dots,\alpha_n)\in\mathbb{N}^n$ and $\partial_x f(x)$ denotes the column vector $(\partial_{x_1}f(x),\dots,\partial_{x_n}f(x))$. If $h\in\mathbb{N}$, then the space of the $h$ times continuously differentiable real-valued functions on $\mathcal{O}$ is denoted by $C^{m}(\mathcal{O})$. The subspace of $C^{h}(\mathcal{O})$ of those functions $f$ whose derivatives $\partial_x^\alpha f$ of order $|\alpha|\le h$ can be extended to a continuous function on $\mathrm{cl}\mathcal{O}$ is denoted $C^h(\mathrm{cl}\mathcal{O})$.

 The space of $n\times n$ real matrices is denoted by $M_n(\mathbb{R})$. If $M\in M_n(\mathbb{R})$ then $M^t$ denotes the transpose matrix of $M$.

For the 
definition and properties of analytic operators, we refer to Prodi and 
Ambrosetti~\cite[p. 89]{PA73} and to Deimling \cite[p.~150]{De85}. In the sequel we shall need the following classical lemma.

\begin{lem}\label{fslem1}
Let $h_1,h_2\in\mathbb{N}\setminus\{0\}$. Let $\mathcal{X}\subseteq\mathbb{R}^{h_1}$, $\mathcal{Y}\subseteq\mathbb{R}^{h_2}$. Assume that $\mathcal{Y}$ is compact. Let $\tau$ be a finite measure on the measurable subsets of $\mathcal{Y}$. Let $f$ be a real analytic function from $\mathcal{X}\times\mathcal{Y}$ to $\mathbb{R}$. Let $F$ be the function from $\mathcal{X}$ to $\mathbb{R}$ defined by $F(x)\equiv\int_{\mathcal{Y}}f(x,y)\,d\tau_y$ for all $x\in\mathcal{X}$. Then $F$  is real analytic.
\end{lem}  Here, we understand that a function $f$ defined on subset $\mathcal{X}$ of a Banach space is real analytic if $f$ is the restriction to $\mathcal{X}$ of a real analytic function defined on an open neighborhood of $\mathcal{X}$.

\section{The G\"unter tangential derivative on $\partial\mathbb{B}_n$}\label{sec2.1}
Let $\mathcal{O}$ be an open subset of a Banach space $\mathcal{B}$. Let $g$ be a real analytic function from $\partial\mathbb{B}_n\times\mathcal{O}$ to $\mathbb{R}$. Then, by definition of analyticity there exist an open neighborhood $\mathcal{U}$ of $\partial\mathbb{B}_n$ in $\mathbb{R}^n$ and a real analytic map $G$ from $\mathcal{U}\times\mathcal{O}$ to $\mathbb{R}$ such that $g=G_{|\partial\mathbb{B}^n\times\mathcal{O}}$ (see also Section \ref{Notation}). The $j$-G\"unter tangential derivative $\mathscr{D}_{\theta_j}g(\theta,b)$  of $g$ at $(\theta,b)\in\partial\mathbb{B}_n\times\mathcal{O}$ is defined by 
\begin{equation}\label{sec2eq1}
\mathscr{D}_{\theta_j}g(\theta,b)\equiv(\partial_{x_j}G)(\theta,b)-\theta_j\sum_{l=1}^n\theta_l(\partial_{x_l} G)(\theta,b)
\end{equation} for all $j\in\{1,\dots,n\}$
 (cf., {\it e.g.}, G\"unter \cite{Gun57}, Kupradze {\it et al.} \cite{KGBB79}, Duduchava {\it et al.} \cite{DuMiMi06}).  As is well known,  $\mathscr{D}_{\theta_j}g(\theta,b)$  does not depend on the particular choice of the extension $G$ and of the open neighborhood $\mathcal{U}$ (see G\"unter \cite[Chap.~1]{Gun57}).
Then we denote by $\mathscr{D}_{\theta}g(\theta,b)$ the column vector $(\mathscr{D}_{\theta_1}g(\theta,b),\dots,\mathscr{D}_{\theta_n}g(\theta,b))$ and we define $\mathscr{D}^\alpha_{\theta}g(\theta,b)\equiv\mathscr{D}^{\alpha_1}_{\theta_1}\dots\mathscr{D}^{\alpha_n}_{\theta_n}g(\theta,b)$ for all $\alpha\in\mathbb{N}^n$. 
By definition \eqref{sec2eq1} one deduces the validity of the following lemma.
\begin{lem}\label{sec2lem3} 
Let $\mathcal{O}$ be an open subset of a Banach space $\mathcal{B}$. Let $g$ be a real analytic map from $\partial\mathbb{B}_n\times\mathcal{O}$ to $\mathbb{R}$. 
Let $\alpha\in\mathbb{N}^n$. Then $\mathscr{D}_\theta^\alpha g$ is real analytic from  $\partial\mathbb{B}_n\times\mathcal{O}$ to $\mathbb{R}$.
\end{lem}
\proof
If $\mathcal{U}$ is an open neighborhood of $\partial\mathbb{B}_n$ and $G$ from $\mathcal{U}\times\mathcal{O}$ to $\mathbb{R}$ is real analytic, then the map from $\mathcal{U}\times\mathcal{O}$ to $\mathbb{C}$ which takes $(x,b)$ to 
\[
(\partial_{x_j}G)(x,b)-x_j \sum_{l=1}^n x_l(\partial_{x_l} G)(x,b)
\]  is real analytic for all $j\in\{1,\dots,n\}$. Then the validity of the lemma follows by \eqref{sec2eq1} and by a standard induction argument.\qed  

\medskip

 Let now $f$ be a real analytic function from $\partial\mathbb{B}_n\times]0,+\infty[$ to $\mathbb{R}$. Let $j\in\{1,\dots,n\}$. Then we observe that
\begin{equation}\label{sec2eq2}
\partial_{x_j}\bigl(f(x/|x|,|x|)\bigr)=\frac{1}{|x|}\bigl(\mathscr{D}_{\theta_j}f\bigr)(x/|x|,|x|)+\frac{x_j}{|x|}\bigl(\partial_r f\bigr)(x/|x|,|x|)\quad\forall x\in\mathbb{R}^n\setminus\{0\}\,, 
\end{equation} where $\partial_r f$ denotes the partial derivative of $f$ with respect to the variable in $]0,+\infty[$.

\section{The auxiliary function $v$}\label{vsec}

As a first step in the construction of our real analytic function $S$ as in \eqref{Sanal}, \eqref{Sfs}, we will show in this section the existence and uniqueness of a  real analytic function $v$ from $\mathscr{E}_R(2k,n)\times\mathbb{R}^n\times\partial\mathbb{B}_{n}\times\mathbb{R}$ to $\mathbb{R}$ such that $\mathbf{L}[\mathbf{a}]v(\mathbf{a},x,\xi,t)=1$ for all $(\mathbf{a},x,\xi,t)\in\mathscr{E}_R(2k,n)\times\mathbb{R}^n\times\partial\mathbb{B}_{n}\times\mathbb{R}$ and $\partial_x^\alpha v(\mathbf{a},x,\xi,t)=0$ for all $(\mathbf{a},x,\xi,t)\in\mathscr{E}_R(2k,n)\times\mathbb{R}^n\times\partial\mathbb{B}_{n}\times\mathbb{R}$ with $x\cdot\xi=t$ and all $\alpha\in\mathbb{N}^n$ with $|\alpha|\le 2k-1$. Then we will investigate some properties of such a function $v$. 

We introduce the following notation.  If  $\mathbf{a}\equiv(a_\alpha)_{\alpha\in N(2k,n)}\in\mathscr{E}_R(2k,n)$,  then we set
\[
\begin{split}
&P[\mathbf{a}](\xi)=P[\mathbf{a}](\xi_{1},\dots,\xi_{n})\equiv\sum_{\alpha=(\alpha_1,\dots,\alpha_n)\in N(2k,n)}a_{\alpha}\xi_1^{\alpha_1}\dots\xi_n^{\alpha_n}\,,\\
&P_0[\mathbf{a}](\xi)=P_0[\mathbf{a}](\xi_{1},\dots,\xi_{n})\equiv\sum_{\alpha=(\alpha_1,\dots,\alpha_n)\in \mathbb{N}^n\,,\,|\alpha|=2k}a_{\alpha}\xi_1^{\alpha_1}\dots\xi_n^{\alpha_n}\,.
\end{split}
\] So that,  $P[\mathbf{a}]$ is a real polynomial  in $n$ variables of degree $2k$ and $P_0[\mathbf{a}]$ is the homogeneous term of $P[\mathbf{a}]$ of degree $2k$ (the so-called ``principal term'' of $P[\mathbf{a}]$). Then we have the following lemma.

\begin{lem}\label{v}
Let $\mathbf{a}\equiv(a_\alpha)_{\alpha\in N(2k,n)}\in\mathscr{E}_R(2k,n)$.  Let 
\begin{equation}\label{veq0}
\rho\in\mathbb{R}\quad\text{and}\quad \rho\ge 1+\frac{\sum_{|\alpha|\le 2k-1}|a_\alpha|}{\inf_{\xi\in\partial\mathbb{B}_n}|P_0[\mathbf{a}](\xi)|}\,.
\end{equation}  Let $\mathbb{D}_\rho\equiv\{z\in\mathbb{C}\,:\,|z|<\rho\}$. Let $v_\mathbf{a}$ be the function from $\mathbb{R}^n\times\partial\mathbb{B}_n\times\mathbb{R}$ to $\mathbb{R}$ defined by  
\begin{equation}\label{veq1}
v_\mathbf{a}(x,\xi,t)\equiv\frac{1}{2\pi i}\int_{\partial\mathbb{D}_\rho}\frac{e^{(x\cdot\xi-t)\zeta}}{\zeta P[\mathbf{a}](\zeta\xi)}\ d\zeta\quad\forall (x,\xi,t)\in\mathbb{R}^n\times\partial\mathbb{B}_n\times\mathbb{R}\,, 
\end{equation} where $i$ denotes the imaginary unit, namely $i^2=-1$, and where $\zeta\xi$ denotes the complex vector $(\zeta\xi_1,\dots,\zeta\xi_n)$. Here we understand that the line integral in \eqref{veq1} is taken over the parametrization $\rho e^{is}$, $s\in[0,2\pi[$.

 Then $v_\mathbf{a}$ is the unique real analytic function from $\mathbb{R}^n\times\partial\mathbb{B}_n\times\mathbb{R}$ to $\mathbb{R}$ such that 
$\mathbf{L}[\mathbf{a}]v_\mathbf{a}(x,\xi,t)=\sum_{|\alpha|\le 2k}a_\alpha\partial_x^\alpha v_\mathbf{a}(x,\xi,t)=1$ for all $(x,\xi,t)\in\mathbb{R}^n\times\partial\mathbb{B}_n\times\mathbb{R}$ and 
$\partial_x^{\alpha}v_\mathbf{a}(x,\xi,t)=0$  for all $(x,\xi,t)\in\mathbb{R}^n\times\partial\mathbb{B}_n\times\mathbb{R}$ with $x\cdot\xi=t$ and all $\alpha\in\mathbb{N}^n$ with $|\alpha|\le 2k-1$.
\end{lem}
\proof By the membership of $\mathbf{a}$ in $\mathscr{E}_R(2k,n)$ one deduces that 
\[
\inf_{\xi\in\partial\mathbb{B}_n}|P_0[\mathbf{a}](\xi)|>0\,.
\] Now let $\zeta\in\mathbb{C}\setminus\{0\}$ and $\xi\in\partial\mathbb{B}_n$. Assume that $P[\mathbf{a}](\zeta\xi)=0$. Then $P_0[\mathbf{a}](\xi)\zeta^{2k}=-(P[\mathbf{a}](\zeta\xi)-P_0[\mathbf{a}](\zeta\xi))$ and thus 
$\zeta=-{\sum_{|\alpha|\le 2k-1}a_\alpha\xi^\alpha\zeta^{|\alpha|-(2k-1)}}/{P_0[\mathbf{a}](\xi)}$. It follows that either $|\zeta|<1$ or $|\zeta|\le{\sum_{|\alpha|\le 2k-1}|a_\alpha|}/{|P_0[\mathbf{a}](\xi)|}$. Thus condition \eqref{veq0} implies that $|\zeta|<\rho$ and one concludes that the polynomial $P[\mathbf{a}](\zeta\xi)$ has no complex zeros $\zeta$ outside of the open disk $\mathbb{D}\rho$, for all $\xi\in\partial\mathbb{B}_n$.   Then Lemma \ref{fslem1} and standard calculus in Banach space imply that the function $v_\mathbf{a}$ defined by \eqref{veq1} is real analytic from $\mathbb{R}^n\times\partial\mathbb{B}_n\times\mathbb{R}$ to $\mathbb{R}$. A straightforward calculation shows that $\mathbf{L}[\mathbf{a}]e^{(x\cdot\xi-t)\zeta}=\sum_{|\alpha|\le 2k}a_\alpha\partial_x^\alpha e^{(x\cdot\xi-t)\zeta}=P[\mathbf{a}](\zeta\xi)$. Then, by standard theorems on differentiation under integral sign and by Cauchy integral formula one has
\begin{equation}\label{veq2}
\mathbf{L}[\mathbf{a}]v_\mathbf{a}(x,\xi,t)=\frac{1}{2\pi i}\int_{\partial\mathbb{D}_\rho}\frac{\mathbf{L}[\mathbf{a}]\left(e^{(x\cdot\xi-t)\zeta}\right)}{\zeta P[\mathbf{a}](\zeta\xi)}\ d\zeta=\frac{1}{2\pi i}\int_{\partial\mathbb{D}_\rho}\frac{1}{\zeta}\ d\zeta=1
\end{equation} for all $(x,\xi,t)\in\mathbb{R}^n\times\partial\mathbb{B}_n\times\mathbb{R}$. By the equality $e^{(x\cdot\xi-t)\zeta}=\sum_{j=0}^{\infty}{(x\cdot\xi-t)^j\zeta^j}/{j!}$ and by standard theorems on summation under integral sign one has
\begin{equation}\label{veq3}
v_\mathbf{a}(x,\xi,t)=\sum_{j=0}^\infty\frac{a_j(\xi)}{j!}(x\cdot\xi-t)^j\quad\forall (x,\xi,t)\in\mathbb{R}^n\times\partial\mathbb{B}_n\times\mathbb{R}
\end{equation} with
\begin{equation}\label{veq4}
a_j(\xi)\equiv\frac{1}{2\pi i}\int_{\partial\mathbb{D}_\rho}\frac{\zeta^{j-1}}{P[\mathbf{a}](\zeta\xi)}\ d\zeta\quad\forall \xi\in\partial\mathbb{B}_n\,,\,j\in\mathbb{N}\,.
\end{equation} 
Now let $\xi\in\partial\mathbb{B}_n$. Let  $g_{\xi}(\zeta)\equiv 1/{P[\mathbf{a}](\xi/\zeta)}$ for all $\zeta\in\mathbb{C}$. Then $g_{\xi}$ is holomorphic in a open neighborhood of $\{z\in\mathbb{C}\,:\,|z|\le1/\rho\}$ and by Cauchy integral formula one verifies that $a_j(\xi)=(1/j!)(\partial_{\zeta}^jg_{\xi})(0)$ for all $j\in\mathbb{N}$ (cf.~equality \eqref{veq4}). One deduces that 
\begin{equation}\label{111201_fs_eqn3}
a_j(\xi)\in\mathbb{R}\qquad\forall\xi\in\partial\mathbb{B}_ n\,,\,j\in\mathbb{N}
\end{equation} and that 
\begin{equation}\label{veq5}
a_j=0\quad \forall j\in\{0,1,\dots,2k-1\}\,,\qquad a_{2k}(\xi)=1/P_{0}[\mathbf{a}](\xi)\quad\forall \xi\in\partial\mathbb{B}_n\,.
\end{equation} Then, by the equalities in \eqref{veq2}, \eqref{veq3}, and \eqref{veq5}, and by standard calculus in Banach space one verifies that $\mathbf{L}[\mathbf{a}]v_\mathbf{a}(x,\xi,t)=1$ for all $(x,\xi,t)\in\mathbb{R}^n\times\partial\mathbb{B}_n\times\mathbb{R}$ and that  $\partial_x^{\alpha}v_\mathbf{a}(x,\xi,t)=0$ for all $(x,\xi,t)\in\mathbb{R}^n\times\partial\mathbb{B}_n\times\mathbb{R}$ with $x\cdot\xi=t$ and all $\alpha\in\mathbb{N}^n$ with $|\alpha|\le 2k-1$. The uniqueness of the function $v_\mathbf{a}$ is an immediate consequence of Cauchy--Kovalevskaya Theorem. The statement is now proved. 
\qed

\medskip

We are now ready to show in the following Proposition \ref{va} the existence and uniqueness of the auxiliary function $v$.

\begin{prop}\label{va} There exist a unique real analytic function $v$ from $\mathscr{E}_R(2k,n)\times\mathbb{R}^n\times\partial\mathbb{B}_{n}\times\mathbb{R}$ to $\mathbb{R}$ such that $\mathbf{L}[\mathbf{a}]v(\mathbf{a},x,\xi,t)= P[\mathbf{a}](\partial_{x_{1}},\dots,\partial_{x_{n}})v(\mathbf{a},x,\xi,t)=1$ for all $(\mathbf{a},x,\xi,t)\in\mathscr{E}_R(2k,n)\times\mathbb{R}^n\times\partial\mathbb{B}_{n}\times\mathbb{R}$ and $\partial_x^\alpha v(\mathbf{a},x,\xi,t)=0$ for all $(\mathbf{a},x,\xi,t)\in\mathscr{E}_R(2k,n)\times\mathbb{R}^n\times\partial\mathbb{B}_{n}\times\mathbb{R}$ with $x\cdot\xi=t$ and for all $\alpha\in\mathbb{N}^n$ with $|\alpha|\le 2k-1$. 
\end{prop}
\proof
We set
\[
\mathscr{E}_{R,l}(2k,n)\equiv\biggl\{\mathbf{a}\in\mathscr{E}_R(2k,n)\,:\,\sum_{|\alpha|\le 2k-1}|a_\alpha|<l\quad\text{and  }\inf_{\xi\in\partial\mathbb{B}_n}|P_0[\mathbf{a}](\xi)|>(1/l)\biggr\}
\] for all $l\in\mathbb{N}\setminus\{0\}$. Then it is easily verified that $\mathscr{E}_{R,l}(2k,n)$ is an open subset of $\mathscr{E}_R(2k,n)$, that  $\cup_{l=1}^\infty\mathscr{E}_{R,l}(2k,n)=\mathscr{E}_{R}(2k,n)$, and that $\mathscr{E}_{R,l}(2k,n)\subseteq \mathscr{E}_{R,l+1}(2k,n)$ for all $l\in\mathbb{N}\setminus\{0\}$.
Now let $l\in\mathbb{N}\setminus\{0\}$ be fixed. Let $\rho_l\equiv 1+l^2$ and $\mathbb{D}_{\rho_l}\equiv\{z\in\mathbb{C}\,:\,|z|<\rho_l\}$.  Let $v_l$ denote the function from $\mathscr{E}_{R,l}(2k,n)\times\mathbb{R}^n\times\partial\mathbb{B}_n\times\mathbb{R}$ to $\mathbb{R}$ defined by
\begin{equation}\label{fseqn0.5}
\begin{split}
&v_l(\mathbf{a},x,\xi,t)\\
&\ \equiv\frac{1}{2\pi i}\int_{\partial\mathbb{D}_{\rho_l}}\frac{e^{(x\cdot\xi-t)\zeta}}{\zeta P[\mathbf{a}](\zeta\xi)}\ d\zeta\quad\forall(\mathbf{a},x,\xi,t)\in\mathscr{E}_{R,l}(2k,n)\times\mathbb{R}^n\times\partial\mathbb{B}_n\times\mathbb{R}\,.
\end{split}
\end{equation} 
One verifies that 
\[
\rho_l\ge 1+\sum_{|\alpha|\le 2k-1}|a_\alpha|/\inf_{\xi\in\partial\mathbb{B}_n}P_0[\mathbf{a}](\xi)\qquad\forall \mathbf{a}\in\mathscr{E}_{R,l}(2k,n)\,.
\] Then $P[\mathbf{a}](\zeta\xi)\neq 0$ for all $\zeta\in\mathbb{C}$ with $|\zeta|=\rho_l$ and all $(\mathbf{a},\xi)\in \mathscr{E}_{R,l}(2k,n)\times\partial\mathbb{B}_n$  (see also the proof of Lemma \ref{v}). Thus, Lemma \ref{fslem1} implies that $v_l$ is real analytic from $\mathscr{E}_{R,l}(2k,n)\times\mathbb{R}^n\times\partial\mathbb{B}_n\times\mathbb{R}$ to $\mathbb{R}$. Moreover, Lemma \ref{v} implies that $v_l$ is the unique real analytic function  from $\mathscr{E}_{R,l}(2k,n)\times\mathbb{R}^n\times\partial\mathbb{B}_n\times\mathbb{R}$ to $\mathbb{R}$ such that $\mathbf{L}[\mathbf{a}]v_l(\mathbf{a},x,\xi,t)=1$ for all $(\mathbf{a},x,\xi,t)\in \mathscr{E}_{R,l}(2k,n)\times\mathbb{R}^n\times\partial\mathbb{B}_n\times\mathbb{R}$ and such that $\partial_x^\alpha v_l(\mathbf{a},x,\xi,t)=0$ for all $(\mathbf{a},x,\xi,t)\in \mathscr{E}_{R,l}(2k,n)\times\mathbb{R}^n\times\partial\mathbb{B}_n\times\mathbb{R}$ with $x\cdot\xi=t$ and for all $\alpha\in\mathbb{N}^n$ with $|\alpha|\le 2k-1$. Hence $v_l(\mathbf{a},x,\xi,t)=v_{l'}(\mathbf{a},x,\xi,t)$ for all $(\mathbf{a},x,\xi,t)\in\mathscr{E}_{R,l}(2k,n)\times\mathbb{R}^n\times\partial\mathbb{B}_n\times\mathbb{R}$ and  all $l,l'\in\mathbb{N}\setminus\{0\}$ such that $l\le l'$. Thus we can define the function $v$ from $\mathscr{E}_R(2k,n)\times\mathbb{R}^n\times\partial\mathbb{B}_n\times\mathbb{R}$ to $\mathbb{R}$ by setting $v(\mathbf{a},x,\xi,t)\equiv v_{l}(\mathbf{a},x,\xi,t)$ for all $(\mathbf{a},x,\xi,t)\in\mathscr{E}_{R,l}(2k,n)\times\mathbb{R}^n\times\partial\mathbb{B}_n\times\mathbb{R}$ and all $l\in\mathbb{N}\setminus\{0\}$. Then $v$ satisfies the conditions in the statement of the proposition.\qed

\medskip

In the following Propositions \ref{wa} and  \ref{aja} we investigate some further properties of the auxiliary function $v$.

\begin{prop}\label{wa} Let $v$ be the function from $\mathscr{E}_R(2k,n)\times\mathbb{R}^n\times\partial\mathbb{B}_{n}\times\mathbb{R}$ to $\mathbb{R}$ of Proposition \ref{va}. 
Then there exists a unique real analytic function $w$ from $\mathscr{E}_R(2k,n)\times\partial\mathbb{B}_n\times\mathbb{R}$ to $\mathbb{R}$ such that $v(\mathbf{a},x,\xi,t)=(x\cdot\xi-t)^{2k}w(\mathbf{a},\xi,x\cdot\xi-t)$ for all $(\mathbf{a},x,\xi,t)\in \mathscr{E}_R(2k,n)\times\mathbb{R}^n\times\partial\mathbb{B}_{n}\times\mathbb{R}$.
\end{prop}
\proof
Let  $\tilde{v}$ be the function from $\mathscr{E}_R(2k,n)\times\partial\mathbb{B}_n\times\mathbb{R}$ to $\mathbb{R}$ which takes $(\mathbf{a},\xi,t)$ to $\tilde{v}(\mathbf{a},\xi,t)\equiv v(\mathbf{a},0,\xi,-t)$. Then $\tilde{v}$ is real analytic from $\mathscr{E}_R(2k,n)\times\partial\mathbb{B}_n\times\mathbb{R}$ to $\mathbb{R}$ and a straightforward verification shows that
\[
v(\mathbf{a},x,\xi,t)=\tilde{v}(\mathbf{a},\xi,x\cdot\xi-t)\qquad\forall(\mathbf{a},x,\xi,t)\in\mathscr{E}_R(2k,n)\times\mathbb{R}^n\times\partial\mathbb{B}_n\times\mathbb{R}
\] (see also \eqref{fseqn0.5}). Hence $\partial_t^j\tilde{v}(\mathbf{a},\xi,0)=0$ for all $(\mathbf{a},\xi)\in\mathscr{E}_R(2k,n)\times\partial\mathbb{B}_n$ and   $j\in\{0,\dots,2k-1\}$.  Then, by standard properties of real analytic functions one can prove that  there exists a unique real analytic function $w$ from $\mathscr{E}_R(2k,n)\times\partial\mathbb{B}_n\times\mathbb{R}$ to $\mathbb{R}$ such that $\tilde{v}(\mathbf{a},\xi,t)=t^{2k}w(\mathbf{a},\xi,t)$ for all $(\mathbf{a},\xi,t)\in\mathscr{E}_R(2k,n)\times\partial\mathbb{B}_n\times\mathbb{R}$. The function $w$ satisfies the conditions in the statement.
\qed

\medskip

\begin{prop}\label{aja}  Let $v$ be the function from $\mathscr{E}_R(2k,n)\times\mathbb{R}^n\times\partial\mathbb{B}_{n}\times\mathbb{R}$ to $\mathbb{R}$ of Proposition \ref{va}. Then there exists a sequence $\{a_j\}_{j\in\mathbb{N}}$ of real analytic functions from $\mathscr{E}_R(2k,n)\times\partial\mathbb{B}_{n}$ to $\mathbb{R}$ such that 
\[
v(\mathbf{a},x,\xi,t)=\sum_{j=0}^\infty \frac{a_j(\mathbf{a},\xi)}{j!}(x\cdot\xi-t)^j\quad\forall(\mathbf{a},x,\xi,t)\in\mathscr{E}_R(2k,n)\times\mathbb{R}^n\times\partial\mathbb{B}_n\times\mathbb{R}\,,
\] where the series converges absolutely and uniformly in the compact subsets of $\mathscr{E}_R(2k,n)\times\mathbb{R}^n\times\partial\mathbb{B}_n\times\mathbb{R}$. Moreover, $a_j=0$ for $j\le 2k-1$ and $a_{2k}(\mathbf{a},\xi)=1/P_0[\mathbf{a}](\xi)$ for all $(\mathbf{a},\xi)\in\mathscr{E}_R(2k,n)\times\partial\mathbb{B}_{n}$.
\end{prop}
\proof Let $\mathscr{E}_{R,l}(2k,n)$, $\rho_l$, $\mathbb{D}_{\rho_l}$, $v_l$ be defined as in the proof of Proposition \ref{va} for all $l\in\mathbb{N}\setminus\{0\}$. Let
\begin{equation}\label{111201_fs_eq2}
a_{j,l}(\mathbf{a},\xi)\equiv\frac{1}{2\pi i}\int_{\partial\mathbb{D}_{\rho_l}}\frac{\zeta^{j-1}}{P[\mathbf{a}](\zeta\xi)}\ d\zeta\quad\forall (\mathbf{a},\xi)\in\mathscr{E}_{R,l}(2k,n)\times\partial\mathbb{B}_n\,,\,j\in\mathbb{N}\,.
\end{equation}  Then $a_{j,l}(\mathbf{a},\xi)\in\mathbb{R}$ for all $(\mathbf{a},\xi)\in\mathscr{E}_{R,l}(2k,n)\times\partial\mathbb{B}_n$ and $j,l\in\mathbb{N}$, $l\ge 1$ (cf.~\eqref{veq4} and \eqref{111201_fs_eqn3}).
Moreover, by arguing so as in the proof of Proposition \ref{va} one verifies that the functions $a_{j,l}$ are real analytic from $\mathscr{E}_{R,l}(2k,n)\times\partial\mathbb{B}_n$ to $\mathbb{R}$ for all $j,l\in\mathbb{N}$, $l\ge 1$.  By the membership of $\mathbf{a}$ in $\mathscr{E}_{R,l}(2k,n)$ one proves that $|P[\mathbf{a}](\zeta\xi)|\ge (1+l^2)^{2k-1}l^{-1}$ for all $\zeta\in\partial\mathbb{D}_{\rho_l}$, $\xi\in\partial\mathbb{B}_n$. Then, by definition \eqref{111201_fs_eq2} and by a straightforward calculation one verifies that 
\begin{equation}\label{ajless}
|a_{j,l}(\mathbf{a},\xi)|\le l(1+l^2)^{j+1-2k}\qquad\forall (\mathbf{a},\xi)\in\mathscr{E}_{R,l}(2k,n)\times\partial\mathbb{B}_n\,.
\end{equation}
Thus, by \eqref{fseqn0.5}, by the inequality in \eqref{ajless}, and by a straightforward calculation one can prove that 
\[
v_l(\mathbf{a},x,\xi,t)=\sum_{j=0}^\infty \frac{a_{j,l}(\mathbf{a},\xi)}{j!}(x\cdot\xi-t)^j\quad\forall(\mathbf{a},x,\xi,t)\in\mathscr{E}_{R,l}(2k,n)\times\mathbb{R}^n\times\partial\mathbb{B}_n\times\mathbb{R},
\] for all $l\in\mathbb{N}\setminus\{0\}$, where the series converges absolutely and uniformly in the compact subsets of  $\mathscr{E}_{R,l}(2k,n)\times\mathbb{R}^n\times\partial\mathbb{B}_n\times\mathbb{R}$. Since  $v_l(\mathbf{a},x,\xi,t)=v_{l'}(\mathbf{a},x,\xi,t)$ for all $(\mathbf{a},x,\xi,t)\in\mathscr{E}_{R,l}(2k,n)\times\mathbb{R}^n\times\partial\mathbb{B}_n\times\mathbb{R}$ and  all $l,l'\in\mathbb{N}\setminus\{0\}$ with $l\le l'$, one deduces that $a_{j,l}(\mathbf{a},\xi)=a_{j,l'}(\mathbf{a},\xi)$ for all $(\mathbf{a},\xi)\in\mathscr{E}_{R,l}(2k,n)\times\partial\mathbb{B}_n$ and  all $j,l,l'\in\mathbb{N}$ with $1\le l\le l'$ (cf.~proof of Proposition \ref{va}). Hence, one can define the function $a_j$ from $\mathscr{E}_{R}(2k,n)\times\partial\mathbb{B}_n$ to $\mathbb{R}$ by setting $a_j(\mathbf{a},\xi)\equiv a_{j,l}(\mathbf{a},\xi)$ for all $(\mathbf{a},\xi)\in\mathscr{E}_{R,l}(2k,n)\times\partial\mathbb{B}_n$ and all $l\in\mathbb{N}\setminus\{0\}$. Then the sequence of functions $\{a_j\}_{j\in\mathbb{N}}$ satisfies the conditions in the statement of the proposition (see also \eqref{veq5}).\qed\medskip

\section{The auxiliary functions $W_0$, $W_1$, and $W_2$}\label{Wsec}

We denote by $W_0$, $W_1$, and $W_2$ the functions from $\mathscr{E}_R(2k,n)\times(\mathbb{R}^n\setminus\{0\})$ to $\mathbb{R}$ defined by 
\begin{eqnarray}
\label{Seq1}
&&W_0(\mathbf{a},x)\equiv\frac{1}{4(2\pi i)^{n-p_n}}\int_{\partial\mathbb{B}_n}\int_0^{x\cdot\xi}v(\mathbf{a},x,\xi,t)\,\mathrm{sgn}\,t\ dt\,d\sigma_\xi\,,\\
\label{Seq2}
&&W_1(\mathbf{a},x)\equiv -\frac{1}{(2\pi i)^{n-p_n}}\int_{\partial\mathbb{B}_n}v(\mathbf{a},x,\xi,0)\log|x\cdot\xi|\ d\sigma_\xi\,,\\
\label{Seq3}
&&W_2(\mathbf{a},x)\equiv -\frac{1}{(2\pi i)^{n-p_n}}\int_{\partial\mathbb{B}_n} \int_0^{x\cdot\xi}\frac{v(\mathbf{a},x,\xi,t)-v(\mathbf{a},x,\xi,0)}{t}\ dt\,d\sigma_\xi\,,\quad
\end{eqnarray} for all $x\in\mathbb{R}^n\setminus\{0\}$ and all $\mathbf{a}\in\mathscr{E}_R(2k,n)$, where   $p_n\equiv 1$ if $n$ is odd and $p_n\equiv 0$ if $n$ is even and where $v$ is the function of Proposition \ref{va}.

In the following Proposition \ref{S} we construct by means of the functions $W_0$, $W_1$, and $W_2$  a particular fundamental solution $S_\mathbf{a}$ for the partial differential  operator $\mathbf{L}[\mathbf{a}]$, $\mathbf{a}\in\mathscr{E}_R(2k,n)$. The validity of Proposition \ref{S} follows by the results of John in \cite[Chap.~3]{J55}. For the sake of completeness we include here a proof. 

\begin{prop}\label{S}
 Let $S_\mathbf{a}$ be the distribution on $\mathbb{R}^n$ defined by
\begin{equation}\label{Sad}
S_\mathbf{a}\equiv p_n \Delta^{(n+1)/2}W_0(\mathbf{a},\cdot) +(1-p_n)\Delta^{n/2}(W_1(\mathbf{a},\cdot)+W_2(\mathbf{a},\cdot)) \,.
\end{equation} Then $S_\mathbf{a}$ is a fundamental solution of the operator $\mathbf{L}[\mathbf{a}]$.
\end{prop}
\proof  Let $\mathbf{a}\in\mathscr{E}_R(2k,n)$ be fixed. Note that the functions $W_0(\mathbf{a},\cdot)$, $W_1(\mathbf{a},\cdot)$, and $W_2(\mathbf{a},\cdot)$ extend to continuous functions on the whole of $\mathbb{R}^n$. Accordingly the expression in the right hand side of \eqref{Sad} defines distribution on $\mathbb{R}^n$.

  Now assume that the dimension $n$  is odd, so that $p_n=1$.  By standard theorems on differentiation under integral sign one verifies that 
\[
\partial_x^\alpha W_0(\mathbf{a},x)=\frac{1}{4(2\pi i)^{n-1}}\int_{\partial\mathbb{B}_n}\int_0^{x\cdot\xi}\partial_x^\alpha v(\mathbf{a},x,\xi,t)\,\mathrm{sgn}\,t\ dt\,d\sigma_\xi\quad\forall x\in\mathbb{R}^n\setminus\{0\}\,,
\] for all $\alpha\in\mathbb{N}^n$ with $|\alpha|\le 2k$ (see also Proposition \ref{va}).  Note that the function $\partial_x^\alpha W_0(\mathbf{a},\cdot)$ extends to a continuous function on the whole of $\mathbb{R}^n$  for all $\alpha\in\mathbb{N}^n$ with $|\alpha|\le 2k$. Moreover, one has
\begin{equation}\label{LW0}
\mathbf{L}[\mathbf{a}] W_0(\mathbf{a},x)=\frac{1}{4(2\pi i)^{n-1}}\int_{\partial\mathbb{B}_n}|x\cdot\xi|\,d\sigma_\xi=\frac{(-1)^{(n-1)/2}}{2^n\pi^{(n-1)/2}((n-1)/2)!}|x| 
\end{equation} for all $x\in\mathbb{R}^n\setminus\{0\}$  (cf.~John \cite[p.~9]{J55}). Thus $\mathbf{L}[\mathbf{a}] W_0(\mathbf{a},\cdot)$ and the function in the right hand side of \eqref{LW0} define the same distribution on $\mathbb{R}^n$. Hence, the distribution $\mathbf{L}[\mathbf{a}] W_0(\mathbf{a},\cdot)$ is a fundamental solution of $\Delta^{(n+1)/2}$ (see, {\it e.g.}, John \cite[p.~44]{J55}). It follows that 
\[
\mathbf{L}[\mathbf{a}](\Delta^{(n+1)/2} W_0(\mathbf{a},\cdot))=\Delta^{(n+1)/{2}}(\mathbf{L}[\mathbf{a}] W_0(\mathbf{a},\cdot))=\delta_0
\] in the sense of distributions in $\mathbb{R}^n$.   Thus $\Delta^{(n+1)/2} W_0(\mathbf{a},\cdot)$ is a fundamental solution of $\mathbf{L}[\mathbf{a}]$ and the proposition is verified for $n$ odd.\par
Now assume that the dimension $n$ is even and $p_n=0$. Denote by $W_3$ the function from $\mathbb{R}^n\setminus\{0\}$ to $\mathbb{R}$ defined by 
\[
W_3(\mathbf{a},x)\equiv-\frac{1}{(2\pi i)^{n}}\int_{\partial\mathbb{B}_n}\int_0^{x\cdot\xi}v(\mathbf{a},x,\xi,t)\,t\log|t|\ dt\,d\sigma_\xi\quad\forall x\in\mathbb{R}^n \setminus\{0\}\,.
\] Then by the classical theorem on differentiation under integral sign one verifies that 
\begin{equation}\label{Seq6}
\partial_x^\alpha W_3(\mathbf{a},x)\equiv-\frac{1}{(2\pi i)^{n}}\int_{\partial\mathbb{B}_n}\int_0^{x\cdot\xi}\partial_x^\alpha v(\mathbf{a},x,\xi,t)\,t\log|t|\ dt\,d\sigma_\xi\quad\forall x\in\mathbb{R}^n \setminus\{0\}\,
\end{equation} for all $\alpha\in\mathbb{N}^n$ with $|\alpha|\le 2k$ (see also Proposition \ref{va}). Note that the function $\partial_x^\alpha W_3(\mathbf{a},\cdot)$ extends to a continuous function on the whole of $\mathbb{R}^n$  for all $\alpha\in\mathbb{N}^n$ with $|\alpha|\le 2k$. Moreover, a straightforward calculation shows that  
\begin{equation}\label{LW3}
\begin{split}
&\mathbf{L}[\mathbf{a}] W_3(\mathbf{a},x)\\
&\quad=-\frac{1}{2(2\pi i)^{n-1}}\int_{\partial\mathbb{B}_n}(x\cdot\xi)^2\log|x\cdot\xi|\,d\sigma_\xi+\frac{1}{4(2\pi i)^{n-1}}\int_{\partial\mathbb{B}_n}(x\cdot\xi)^2\,d\sigma_\xi\\
&\quad=\frac{(-1)^{({n}/{2})-1}}{2^{n+1}\pi^{n/2}(n/2)!}|x|^2\log|x|+c_n|x|^2\qquad\qquad\forall x\in\mathbb{R}^n \setminus\{0\}\,,
\end{split}
\end{equation} where $c_n$ is a real constant which depends only on $n$ (cf.~John \cite[p.~9]{J55}).  Thus $\mathbf{L}[\mathbf{a}] W_3(\mathbf{a},\cdot)$ and the function in the right hand side of \eqref{LW3} define the same distribution on $\mathbb{R}^n$.  Hence the distribution $\mathbf{L}[\mathbf{a}] W_3(\mathbf{a},\cdot)$ is a fundamental solution of $\Delta^{(n/2)+1}$ (see, {\it e.g.}, John \cite[p.~44]{J55}, also note that $\Delta^{(n/2)+1}|x|^2=0$). Moreover, by Proposition \ref{aja} one verifies that
\[
\partial_x^\alpha v(\mathbf{a},x,\xi,t)=(-1)^{|\alpha|}\xi^\alpha\partial_t^{|\alpha|}v(\mathbf{a},x,\xi,t)\qquad\forall (x,\xi,t)\in\mathbb{R}^n\times\partial\mathbb{B}_n\times\mathbb{R}\,,\ \alpha\in\mathbb{N}^n\,.
\] Hence, equality \eqref{Seq6} implies that 
\begin{equation}\label{Seq7}
\Delta W_3(\mathbf{a},x)=-\frac{1}{(2\pi i)^{n}}\int_{\partial\mathbb{B}_n}\int_0^{x\cdot\xi}(\partial_t^2 v(\mathbf{a},x,\xi,t))\,t\log|t|\ dt\,d\sigma_\xi\quad\forall x\in\mathbb{R}^n \setminus\{0\}\,.
\end{equation} Then, by integrating by parts the integral in the right hand side of \eqref{Seq7} one deduces that
\[
\begin{split}
\Delta W_3(\mathbf{a},x)=&-\frac{1}{(2\pi i)^{n}}\int_{\partial\mathbb{B}_n}v(\mathbf{a},x,\xi,0)\log|x\cdot\xi|\ d\sigma_\xi\\
&-\frac{1}{(2\pi i)^{n-p_n}}\int_{\partial\mathbb{B}_n} \int_0^{x\cdot\xi}\frac{v(\mathbf{a},x,\xi,t)-v(\mathbf{a},x,\xi,0)}{t}\ dt\,d\sigma_\xi\\
&-\frac{1}{(2\pi i)^{n}}\int_{\partial\mathbb{B}_n}v(\mathbf{a},x,\xi,0)\ d\sigma_\xi\qquad\forall x\in\mathbb{R}^n \setminus\{0\}\,.
\end{split}
\] So that
\[
W_1(\mathbf{a},x)+W_2(\mathbf{a},x)=\Delta W_3(\mathbf{a},x)+\frac{1}{(2\pi i)^{n}}\int_{\partial\mathbb{B}_n}v(\mathbf{a},x,\xi,0)\ d\sigma_\xi\qquad\forall x\in\mathbb{R}^n \setminus\{0\}\,.
\] We now observe that 
\[
\mathbf{L}[\mathbf{a}]\left[\frac{1}{(2\pi i)^{n}}\int_{\partial\mathbb{B}_n}v(\mathbf{a},x,\xi,0)\ d\sigma_\xi\right]=\frac{s_n}{(2\pi i)^{n}}\qquad\forall x\in\mathbb{R}^n\,,
\] where $s_n$ denotes the $n-1$ dimensional measure of $\partial\mathbb{B}_n$. Hence, $\mathbf{L}[\mathbf{a}](W_1(\mathbf{a},\cdot)+W_2(\mathbf{a},\cdot))$ differs from $ \mathbf{L}[\mathbf{a}](\Delta W_3(\mathbf{a},\cdot))$ by a constant function. It follows that  
\[
\begin{split}
&\mathbf{L}[\mathbf{a}]\Delta^{n/2}(W_1(\mathbf{a},\cdot)+W_2(\mathbf{a},\cdot))\\
&\qquad=\mathbf{L}[\mathbf{a}]\Delta^{(n/2)+1}W_3(\mathbf{a},\cdot)=\Delta^{(n/2)+1}\mathbf{L}[\mathbf{a}]W_3(\mathbf{a},\cdot)=\delta_0
\end{split}
\]  in the sense of distributions in $\mathbb{R}^n$. Thus $\Delta^{{n}/{2}} (W_1(\mathbf{a},\cdot)+W_2(\mathbf{a},\cdot))$ is a fundamental solution of $\mathbf{L}[\mathbf{a}]$ and the proposition is verified also for $n$ even. \qed

\medskip

In the following Proposition \ref{ABC} we investigate some further properties of the functions $W_0$, $W_1$, and $W_2$.

\begin{prop}\label{ABC}  There exist real analytic functions $A^\flat$ and $A^\sharp$ from $\mathscr{E}_R(2k,n)\times\partial\mathbb{B}_{n}\times\mathbb{R}$ to $\mathbb{R}$,  and real analytic functions $B^\sharp$ and $C^\sharp$ from $\mathscr{E}_R(2k,n)\times\mathbb{R}^n$ to $\mathbb{R}$ such that 
\[
\begin{split}
&W_0(\mathbf{a},x)=|x|^{2k+1}A^\flat(\mathbf{a},x/|x|,|x|)\,,\\
&W_1(\mathbf{a},x)=|x|^{2k}A^\sharp(\mathbf{a},x/|x|,|x|)+\log|x|\,B^\sharp(\mathbf{a},x)\,,\\
&W_2(\mathbf{a},x)=C^\sharp(\mathbf{a},x)\,,\qquad\qquad\forall (\mathbf{a},x)\in\mathscr{E}_R(2k,n)\times(\mathbb{R}^n\setminus\{0\})\,.
\end{split}
\]
 Moreover,  $\partial^\alpha_x B(\mathbf{a},0)=0$ for all $\mathbf{a}\in\mathscr{E}_R(2k,n)$ and all $\alpha\in\mathbb{N}^n$ with $|\alpha|\le 2k-1$.
\end{prop}
\proof  Let $A^\flat$ be the function from $\mathscr{E}_R(2k,n)\times\partial\mathbb{B}_n\times\mathbb{R}$ to $\mathbb{R}$  defined by 
\begin{equation}\label{fseqn0.6}
{A}^\flat(\mathbf{a},\theta,r)\equiv\frac{1}{4(2\pi i)^{n-p_n}}\int_{\partial\mathbb{B}_n}\int_0^{\theta\cdot\xi}\left(\theta\cdot\xi-s\right)^{2k} w\left(\mathbf{a},\xi,r\left(\theta\cdot\xi-s\right)\right)\mathrm{sgn}\,s\ ds\,d\sigma_\xi
\end{equation} for all $(\mathbf{a},\theta,r)\in\mathscr{E}_R(2k,n)\times\partial\mathbb{B}_n\times\mathbb{R}$, where   $w$ is as in Proposition \ref{wa}. Then,  by definition \eqref{Seq1}, by equality ${v}(\mathbf{a},x,\xi,t)=(x\cdot\xi-t)^{2k}w(\mathbf{a},\xi,x\cdot\xi-t)$, and by a straightforward calculation one verifies that 
\[
W_{0}(\mathbf{a},x)=|x|^{2k+1}A^\flat(\mathbf{a},x/|x|,|x|)\qquad\forall(\mathbf{a},x)\in\mathscr{E}_R(2k,n)\times(\mathbb{R}^n\setminus\{0\})\,.
\] We now show that $A^\flat$ is real analytic from $\mathscr{E}_R(2k,n)\times\partial\mathbb{B}_n\times\mathbb{R}$ to $\mathbb{R}$. By changing variable of integration  in the inner integral in \eqref{fseqn0.6} one obtains 
\begin{equation}\label{ABCeqn1}
\begin{split}
&{A}^\flat(\mathbf{a},\theta,r)\\
&=\frac{1}{4(2\pi i)^{n-p_n}}\int_{\partial\mathbb{B}_n}\int_0^{1}(\theta\cdot\xi)^{2k+1}(1-t)^{2k} w\left(\mathbf{a},\xi,r(\theta\cdot\xi)(1-t)\right)\mathrm{sgn}(\theta\cdot\xi)\ dt\,d\sigma_\xi
\end{split}
\end{equation} Then one introduces a new variable of integration $\xi'$ instead of $\xi$ in \eqref{fseqn0.6}  by a suitable orthogonal substitution. Let $\eta$ be an arbitrary chosen unit vector and consider $\theta$ restricted to the half sphere \
\begin{equation}\label{Bneta+}
\partial\mathbb{B}^+_{n,\eta}\equiv\{\theta\in\partial\mathbb{B}_n\,:\,\theta\cdot\eta>0\}\,.
\end{equation}  Let $T_\eta$ denote the real analytic matrix valued function from $\partial\mathbb{B}^+_{n,\eta}$ to   $M_n(\mathbb{R})$ with $(j,k)$ entry $T_{\eta, jk}$ defined by
\begin{equation}\label{fseqn5.2}
T_{\eta,jk}(\theta)\equiv \delta_{jk}+2\theta_j\eta_k-\frac{(\theta+\eta)_j(\theta+\eta)_k}{1+\theta\cdot\eta}\quad\forall \theta\in\partial\mathbb{B}^+_{n,\eta}\,.
\end{equation}   One verifies that $T_\eta(\theta)$ is an orthogonal matrix and that $T_\eta(\theta)^t\theta=\eta$ for all $\theta\in\partial\mathbb{B}^+_{n,\eta}$. In particular, $\theta\cdot T_\eta(\theta)\xi'=\eta\cdot\xi'$ for all $\xi'\in\mathbb{R}^n$ and all $\partial\mathbb{B}^+_{n,\eta}$.   Then, by taking $\xi=T_\eta(\theta)\xi'$ in \eqref{ABCeqn1} one obtains 
\[
\begin{split}
&A^\flat(\mathbf{a},\theta,r)=\frac{1}{4(2\pi i)^{n-p_n}}\int_{\partial\mathbb{B}_n}\int_0^{1} w\left(\mathbf{a},T_\eta(\theta)\xi',r(\eta\cdot\xi')(1-t)\right)\\
&\qquad\qquad\qquad\qquad\qquad \times(\eta\cdot\xi')^{2k+1}(1-t)^{2k}\mathrm{sgn}(\eta\cdot\xi')\ ds\,d\sigma_{\xi'}
\end{split}
\] for all $(\mathbf{a},\theta,r)\in\mathscr{E}_R(2k,n)\times \partial\mathbb{B}_{n,\eta}^+\times\mathbb{R}$. Observe that the map from $\partial\mathbb{B}_n\times[0,1]$ to $\mathbb{R}$ which takes $(\xi',t)$ to $(\eta\cdot\xi')^{2k+1}(1-t)^{2k}\mathrm{sgn}(\eta\cdot\xi')$ is bounded, so that the measure $(\eta\cdot\xi')^{2k+1}(1-t)^{2k}\mathrm{sgn}(\eta\cdot\xi')dt\,d\sigma_{\xi'}$ on the compact set $\partial\mathbb{B}_n\times[0,1]$ is finite. Moreover, by standard properties of real analytic functions one shows that the map from $\mathscr{E}_R(2k,n)\times \partial\mathbb{B}_{n,\eta}^+\times\mathbb{R}\times\partial\mathbb{B}_n\times[0,1]$ to $\mathbb{R}$ which takes $(\mathbf{a},\theta,r,\xi',t)$ to $w\left(\mathbf{a},T_\eta(\theta)\xi',r(\eta\cdot\xi')(1-t)\right)$ is real analytic. Hence, Lemma \ref{fslem1} implies that the restriction of the function $A^\flat$ to $\mathscr{E}_R(2k,n)\times \partial\mathbb{B}_{n,\eta}^+\times\mathbb{R}$ is real analytic. Since $\eta$ is an arbitrarily chosen unit vector of $\mathbb{R}^n$, it follows that $A^\flat$ is real analytic from $\mathscr{E}_R(2k,n)\times \partial\mathbb{B}_n\times\mathbb{R}$ to $\mathbb{R}$.

Now consider the function $W_1$.  Denote by $A^\sharp$ and $B^\sharp$ the functions from  $\mathscr{E}_R(2k,n)\times \partial\mathbb{B}_n\times\mathbb{R}$ to $\mathbb{R}$ and from $\mathscr{E}_R(2k,n)\times \mathbb{R}^n$ to $\mathbb{R}$, respectively, defined by
\[
A^\sharp(\mathbf{a},\theta,r)\equiv-\frac{1}{(2\pi i)^{n-p_n}}\int_{\partial\mathbb{B}_n}w(\mathbf{a},\xi,r(\theta\cdot\xi))(\theta\cdot\xi)^{2k}\log|\theta\cdot\xi|\,d\sigma_\xi
\] for all $(\mathbf{a},\theta,r)\in \mathscr{E}_R(2k,n)\times \partial\mathbb{B}_n\times\mathbb{R}$,
and
\[
B^\sharp(\mathbf{a},x)\equiv -\frac{1}{(2\pi i)^{n-p_n}}\int_{\partial\mathbb{B}_n}v(\mathbf{a},x,\xi,0)\,d\sigma_\xi\quad\forall (\mathbf{a},\theta,r)\in \mathscr{E}_R(2k,n)\times \mathbb{R}^n\,,
\] where $v$ and $w$ are as in Propositions \ref{va} and \ref{wa}. Then,  by definition \eqref{Seq2}, by equality ${v}(\mathbf{a},x,\xi,t)=(x\cdot\xi-t)^{2k}w(\mathbf{a},\xi,x\cdot\xi-t)$, and by a straightforward calculation one verifies that 
\[
W_1(\mathbf{a},x)=|x|^{2k}A^\sharp(\mathbf{a},x/|x|,|x|)+\log|x|\, B^\sharp(\mathbf{a},x)\ \forall(\mathbf{a},x)\in\mathscr{E}_R(2k,n)\times(\mathbb{R}^n\setminus\{0\})\,.
\] By Lemma \ref{fslem1} one proves that $B^\sharp$ is real analytic from $\mathscr{E}_R(2k,n)\times \mathbb{R}^n$ to $\mathbb{R}$. Then, by equality $\partial_x^\alpha v(\mathbf{a},0,\xi,0)=0$ for all $(\mathbf{a},\xi)\in\mathscr{E}_R(2k,n)\times \partial\mathbb{B}_n$ and all $\alpha\in\mathbb{N}^n$ with $|\alpha|\le 2k-1$ (cf.~Proposition 
\ref{va}), and by  standard theorems on differentiation under integral sign, one shows that $\partial_x^\alpha B^\sharp(\mathbf{a},0)=0$ for all $\mathbf{a}\in\mathscr{E}_R(2k,n)$ and all $\alpha\in\mathbb{N}^n$ with $|\alpha|\le 2k-1$. Now one has to prove that $A^\sharp$ is real analytic from $\mathscr{E}_R(2k,n)\times \partial\mathbb{B}_n\times\mathbb{R}$ to $\mathbb{R}$. To do so, fix a unit vector $\eta\in\partial\mathbb{B}_n$. Then verify that 
\[
A^\sharp(\mathbf{a},\theta,r)=-\frac{1}{(2\pi i)^{n-p_n}}\int_{\partial\mathbb{B}_n}w(\mathbf{a},T_\eta(\theta)\xi',r(\eta\cdot\xi'))(\eta\cdot\xi')^{2k}\log|\eta\cdot\xi'|\,d\sigma_{\xi'}
\] for all $(\mathbf{a},\theta,r)\in \mathscr{E}_R(2k,n)\times \partial\mathbb{B}_{n,\eta}^+\times\mathbb{R}$, where $T_\eta(\theta)$ is the orthogonal matrix introduced in \eqref{fseqn5.2}. The map from $\partial\mathbb{B}_n$ to $\mathbb{R}$ which takes $\xi'$ to $(\eta\cdot\xi')^{2k}\log|\eta\cdot\xi'|$ is bounded, so that the measure $(\eta\cdot\xi')^{2k}\log|\eta\cdot\xi'|\,d\sigma_{\xi'}$ on the compact set $\partial\mathbb{B}_n$ is finite. Moreover, by standard properties of real analytic functions one can verify that the map from $\mathscr{E}_R(2k,n)\times \partial\mathbb{B}_{n,\eta}^+\times\mathbb{R}\times\partial\mathbb{B}_n$ to $\mathbb{R}$ which takes $(\mathbf{a},\theta,r,\xi')$ to $w(\mathbf{a},T_\eta(\theta)\xi',r(\eta\cdot\xi'))$ is real analytic. Hence, Lemma \ref{fslem1} implies that the restriction of the function $A^\sharp$ to $\mathscr{E}_R(2k,n)\times \partial\mathbb{B}_{n,\eta}^+\times\mathbb{R}$ is real analytic. Since $\eta$ is an arbitrarily chosen unit vector of $\mathbb{R}^n$, it follows that $A^\sharp$ is real analytic from $\mathscr{E}_R(2k,n)\times \partial\mathbb{B}_n\times\mathbb{R}$ to $\mathbb{R}$.

Finally consider the function $W_2$. Let $\tilde v$ be the function from $\mathscr{E}_R(2k,n)\times \mathbb{R}^n\times\partial\mathbb{B}_n\times\mathbb{R}$ to $\mathbb{R}$ defined by $\tilde v(\mathbf{a},x,\xi,t)\equiv(v(\mathbf{a},x,\xi,t)-v(\mathbf{a},x,\xi,0))/t$ for all $(\mathbf{a},x,\xi,t)\in \mathscr{E}_R(2k,n)\times \mathbb{R}^n\times\partial\mathbb{B}_n\times\mathbb{R}$. Then, by standard properties of real analytic functions one shows that $\tilde v$ is real analytic from $\mathscr{E}_R(2k,n)\times \mathbb{R}^n\times\partial\mathbb{B}_n\times\mathbb{R}$ to $\mathbb{R}$. Now set
\[
C^\sharp(\mathbf{a},x)\equiv -\frac{1}{(2\pi i)^{n-p_n}}\int_{\partial\mathbb{B}_n}\int_0^1(x\cdot\xi)\tilde v(\mathbf{a},x,\xi,(x\cdot\xi)s)\,ds\, d\sigma_\xi
\] for all $(\mathbf{a},x)\in \mathscr{E}_R(2k,n)\times \mathbb{R}^n$. Then Lemma \ref{fslem1} implies that  $C^\sharp$ is real analytic from $\mathscr{E}_R(2k,n)\times \mathbb{R}^n$ to $\mathbb{R}$ and a straightforward calculation shows that $W_2(\mathbf{a},x)=C^\sharp(\mathbf{a},x)$ for all $(\mathbf{a},x)\in \mathscr{E}_R(2k,n)\times(\mathbb{R}^n\setminus\{0\})$.
\qed\medskip

\section{The functions $S$ and $S_0$}\label{Ssec}

In the following Theorem \ref{Sa} we introduce a real analytic function $S$ which satisfies the conditions in \eqref{Sanal} and \eqref{Sfs}.

\begin{thm}\label{Sa}
Let $S$ be the function from $\mathscr{E}_R(2k,n)\times (\mathbb{R}^n\setminus\{0\})$ to $\mathbb{R}$ defined by
\[
S(\mathbf{a},x)\equiv p_n \Delta^{(n+1)/2}W_0(\mathbf{a},x) +(1-p_n)\Delta^{n/2}(W_1(\mathbf{a},x)+W_2(\mathbf{a},x))
\] for all $(\mathbf{a},x)\in \mathscr{E}_R(2k,n)\times (\mathbb{R}^n\setminus\{0\})$. Then $S$ is real analytic from $\mathscr{E}_R(2k,n)\times (\mathbb{R}^n\setminus\{0\})$ to $\mathbb{R}$, $S(\mathbf{a},\cdot)$ is a locally integrable function on $\mathbb{R}^n$ for all $\mathbf{a}\in \mathscr{E}_R(2k,n)$, and $S(\mathbf{a},\cdot)$  is a fundamental solution of the operator $\mathbf{L}[\mathbf{a}]$ for all $\mathbf{a}\in \mathscr{E}_R(2k,n)$.
\end{thm}
\proof Proposition \ref{ABC} and standard properties of real analytic functions imply that $S$ is real analytic. Then, by Proposition \ref{ABC} and by a straightforward calculation one verifies that $S(\mathbf{a},\cdot)$ is a locally integrable function on $\mathbb{R}^n$ for all $\mathbf{a}\in \mathscr{E}_R(2k,n)$. Hence   $S(\mathbf{a},\cdot)$ defines a distribution on $\mathbb{R}^n$. Finally, by Proposition \ref{ABC} and by a standard argument based on the divergence theorem one verifies that $S(\mathbf{a},\cdot)= p_n \Delta^{(n+1)/2}W_0(\mathbf{a},\cdot) +(1-p_n)\Delta^{n/2}(W_1(\mathbf{a},\cdot)+W_2(\mathbf{a},\cdot))
$ in the sense of distributions on $\mathbb{R}^n$. Hence Proposition \ref{S} implies that $S(\mathbf{a},\cdot)=S_\mathbf{a}$, and thus $S(\mathbf{a},\cdot)$ is a fundamental solution of $\mathbf{L}[\mathbf{a}]$ for all $\mathbf{a}\in \mathscr{E}_R(2k,n)$.\qed

\medskip

In Theorem \ref{fsthm1} here below we provide a detailed expression for $S$.

\begin{thm}\label{fsthm1}  Let $S$ be as in Theorem \ref{Sa}. Then, there exist a real analytic function $A$ from $\mathscr{E}_R(2k,n)\times\partial\mathbb{B}_{n}\times\mathbb{R}$ to $\mathbb{R}$,  and real analytic functions $B$ and $C$ from $\mathscr{E}_R(2k,n)\times\mathbb{R}^n$ to $\mathbb{R}$ such that 
\begin{equation}\label{fsthm1eq0}
S(\mathbf{a},x)= |x|^{2k-n}A(\mathbf{a},x/|x|,|x|)+ \log |x|\,  B(\mathbf{a},x) +C(\mathbf{a},x)
\end{equation}  for all $(\mathbf{a},x)\in\mathscr{E}_R(2k,n)\times(\mathbb{R}^n\setminus\{0\})$.
The functions $B$ and $C$ are identically $0$ if $n$ is odd and there exist    a sequence  $\{f_j\}_{j\in\mathbb{N}}$ of real analytic functions from $\mathscr{E}_R(2k,n)\times\partial\mathbb{B}_n$ to $\mathbb{R}$, and  a family $\{b_\alpha\}_{|\alpha|\ge\sup\{k-n,0\}}$ of real analytic functions from $\mathscr{E}_R(2k,n)$ to $\mathbb{R}$,  such that 
\[
f_{j}(\mathbf{a},-\theta)=(-1)^{j}f_{j}(\mathbf{a},\theta)\quad\forall (\mathbf{a},\theta)\in\mathscr{E}_R(2k,n)\times\partial\mathbb{B}_n\,,
\] and
\begin{eqnarray}\label{thm1eq1}
&&A(\mathbf{a},\theta,r)=\sum_{j=0}^{\infty}r^jf_j(\mathbf{a},\theta)\qquad \forall\ (\mathbf{a},\theta,r)\in\mathscr{E}_R(2k,n)\times\partial\mathbb{B}_{n}\times\mathbb{R},\\
\label{thm1eq2}
&&B(\mathbf{a},x)=\sum_{|\alpha|\ge\sup\{2k-n,0\}}b_{\alpha}(\mathbf{a})x^\alpha
\qquad \forall  (\mathbf{a},x)\in\mathscr{E}_R(2k,n)\times\mathbb{R}^n,
\end{eqnarray}  where the series in \eqref{thm1eq1} and \eqref{thm1eq2} converge absolutely and uniformly in all compact subsets of  $\mathscr{E}_R(2k,n)\times\partial\mathbb{B}_{n}\times\mathbb{R}$ and  of $\mathscr{E}_R(2k,n)\times\mathbb{R}^n$, respectively.
\end{thm}
\proof One proves the theorem separately for $n$ odd and $n$ even. First assume that $n$ is odd. By Proposition \ref{ABC}, one has $W_0(\mathbf{a},x)=|x|^{2k+1}A^\flat(\mathbf{a},x/|x|,|x|)$ with $A^\flat$ real analytic from $\mathscr{E}_R(2k,n)\times\partial\mathbb{B}_n\times\mathbb{R}$ to $\mathbb{R}$. If $h\in\{1,\dots,n\}$, then one verifies that 
\[
\partial_{x_h} W_{0}(\mathbf{a},x)=|x|^{2k}A^{\flat,h}(\mathbf{a},x/|x|,|x|)\qquad\forall(\mathbf{a},x)\in \mathscr{E}_R(2k,n)\times(\mathbb{R}^n\setminus\{0\})\,,
\] where $A^{\flat,h}$ is the function  defined by
\[
A^{\flat,h}(\mathbf{a},\theta,r)\equiv(2k+1)\theta_h A^\flat(\mathbf{a},\theta,r)+\theta_h \,r\partial_{r}A^\flat(\mathbf{a},\theta,r)+\mathscr{D}_{\theta_h}A^\flat(\mathbf{a},\theta,r)
\]  for all $(\mathbf{a},\theta,r)\in \mathscr{E}_R(2k,n)\times\partial\mathbb{B}_{n}\times\mathbb{R}$ (see also \eqref{sec2eq2}). Observe that $A^{\flat,h}$ is real analytic from $\mathscr{E}(2k,n)\times\partial\mathbb{B}_{n}\times\mathbb{R}$ to $\mathbb{R}$ (see also Lemma \ref{sec2lem3}). 
Then, by an induction argument on the order of differentiation one proves that there exists a real analytic map  $A$ from $\mathscr{E}(2k,n)\times\partial\mathbb{B}_{n}\times\mathbb{R}$ to $\mathbb{R}$ such that
\[
\Delta^{(n+1)/2}W_{0}(\mathbf{a},x)=|x|^{2k-n}A(\mathbf{a},x/|x|,|x|)\quad\forall(\mathbf{a},x)\in \mathscr{E}_R(2k,n)\times(\mathbb{R}^n\setminus\{0\})\,.
\]  Hence,  equality \eqref{fsthm1eq0} for $n$ odd holds with $B\equiv 0$ and $C\equiv 0$ (cf.~Theorem \ref{Sa}). To complete the proof for $n$ odd one has  to show the existence of the family of functions $\{f_j\}_{j\in\mathbb{N}}$ as in the statement. By Proposition \ref{aja}, by definition \eqref{Seq1}, by the inequality in \eqref{ajless}, and by the dominated convergence theorem, one has 
\[
W_0(\mathbf{a},x)=\sum_{j=2k}^{\infty}W_{0,j}(\mathbf{a},x)\qquad\forall (\mathbf{a},x)\in\mathscr{E}_R(2k,n)\times(\mathbb{R}^n\setminus\{0\})\,
\] with 
\[
W_{0,j}(\mathbf{a},x)\equiv-\frac{1}{4(2\pi i)^{n-1}}\int_{\partial\mathbb{B}_n}\frac{a_j(\mathbf{a},\xi)}{(j+1)!}\ (x\cdot\xi)^{j+1}\ \mathrm{sgn}(x\cdot\xi)\ d\sigma_\xi
\] for all $(\mathbf{a},x)\in\mathscr{E}_R(2k,n)\times(\mathbb{R}^n\setminus\{0\})$ and $j\in\mathbb{N}$. Let $A^\flat_{j}$ be defined by
\begin{equation}\label{Aflatj}
A^\flat_{j}(\mathbf{a},\theta)\equiv-\frac{1}{4(2\pi i)^{n-1}}\int_{\partial\mathbb{B}_n}\frac{a_j(\mathbf{a},\xi)}{(j+1)!}\ (\theta\cdot\xi)^{j+1}\ \mathrm{sgn}(\theta\cdot\xi)\ d\sigma_\xi
\end{equation} for all $(\mathbf{a},\theta)\in\mathscr{E}_R(2k,n)\times\mathbb{B}_n$ and $j\in\mathbb{N}$. Then 
\[
W_{0,j}(\mathbf{a},x)=|x|^{j+1}A^\flat_j(\mathbf{a},x/|x|)\quad\forall(\mathbf{a},x)\in \mathscr{E}_R(2k,n)\times(\mathbb{R}^n\setminus\{0\})\,,\ j\in\mathbb{N}\,.
\]
We show that $A^\flat_j$ is real analytic. Let $\eta\in\partial\mathbb{B}_n$. Then, one verifies that 
\[
A^\flat_{j}(\mathbf{a},\theta)=-\frac{1}{4(2\pi i)^{n-1}}\int_{\partial\mathbb{B}_n}\frac{a_j(\mathbf{a},T_\eta(\theta)\xi')}{(j+1)!}\ (\eta\cdot\xi')^{j+1}\ \mathrm{sgn}(\eta\cdot\xi')\ d\sigma_{\xi'}
\]   for all $(\mathbf{a},\theta)\in\mathscr{E}_R(2k,n)\times\partial\mathbb{B}^+_{n,\eta}$ and all $j\in\mathbb{N}$, where $T_\eta(\theta)$ is defined as in \eqref{fseqn5.2} (see also \eqref{Bneta+}). Hence, by Lemma \ref{fslem1}, by the boundedness of the function which takes $\xi'\in\partial\mathbb{B}_n$ to $(\eta\cdot\xi')^{j+1}\ \mathrm{sgn}(\eta\cdot\xi')$, and by standard properties of real analytic functions, one proves that the restriction of $A^\flat_{j}$ to  $\mathscr{E}_R(2k,n)\times\partial\mathbb{B}^+_{n,\eta}$ is real analytic. Thus $A^\flat_{j}$ is real analytic from $\mathscr{E}_R(2k,n)\times\partial\mathbb{B}_n$ to $\mathbb{R}$ for all $j\in\mathbb{N}$ (see also the proof of Proposition \ref{ABC} where a similar argument has been exploited to show the real analyticity of the functions $A^\flat$, $A^\sharp$, $B^\sharp$, and $C^\sharp$). Moreover, 
\[
A^\flat_j(\mathbf{a},-\theta)=(-1)^j A^\flat_j(\mathbf{a},\theta)\qquad\forall (\mathbf{a},\theta)\in \mathscr{E}_R(2k,n)\times\partial\mathbb{B}_n\,,\, j\in\mathbb{N}\,.
\]  Then, by inequality \eqref{ajless} and by definition \eqref{Aflatj} one verifies that the series 
\[
\sum_{j=0}^\infty r^{j}A^\flat_{2k+j}(\mathbf{a},\theta)
\] converges absolutely and uniformly in the compact subsets of $\mathscr{E}_R(2k,n)\times\partial\mathbb{B}_n\times\mathbb{R}$. Now let $h\in\{1,\dots,n\}$. Then one has
\[
\partial_{x_h}\bigl( |x|^{2k+1+j}A^\flat_{2k+j}(\mathbf{a},x/|x|)\bigr)= |x|^{2k+j}A^{\flat,h}_{2k+j}(\mathbf{a},x/|x|)
\] for all $(\mathbf{a},x)\in \mathscr{E}_R(2k,n)\times(\mathbb{R}^n\setminus\{0\})$ and $j\in\mathbb{N}$, where  $A^{\flat,h}_{2k+j}$ is the real analytic function from $\mathscr{E}_R(2k,n)\times\partial\mathbb{B}_n\times\mathbb{R}$ to $\mathbb{R}$ defined by
\[
A^{\flat,h}_{2k+j}(\mathbf{a},\theta)\equiv \mathscr{D}_{\theta_h}A^\flat_{2k+j}(\mathbf{a},\theta)+(2k+1+j)\theta_h A^\flat_{2k+j}(\mathbf{a},\theta) 
\] for all $(\mathbf{a},\theta)\in \mathscr{E}_R(2k,n)\times\partial\mathbb{B}_n$ and $j\in\mathbb{N}$ (see also \eqref{sec2eq2} and Lemma \ref{sec2lem3}). Note that 
\[
A^{\flat,h}_j(\mathbf{a},-\theta)=(-1)^{j+1} A^{\flat,h}_j(\mathbf{a},\theta)\qquad\forall (\mathbf{a},\theta)\in \mathscr{E}_R(2k,n)\times\partial\mathbb{B}_n\,,\, j\in\mathbb{N}\,.
\] Moreover, by standard properties of real analytic functions one verifies that the series 
\[
\sum_{j=0}^\infty r^{2k+1+j}A^{\flat,h}_{2k+j}(\mathbf{a},\theta)=(\mathscr{D}_{\theta_h}+r\theta_h\partial_r)\sum_{j=0}^\infty r^{2k+1+j}A^\flat_{2k+j}(\mathbf{a},\theta)
\]  converges absolutely and uniformly in the compact subsets of $\mathscr{E}_R(2k,n)\times\partial\mathbb{B}_n\times\mathbb{R}$, which in turn implies that the series 
\[
\sum_{j=0}^\infty r^{j}A^{\flat,h}_{2k+j}(\mathbf{a},\theta)
\] converges absolutely and uniformly in the compact subsets of $\mathscr{E}_R(2k,n)\times\partial\mathbb{B}_n\times\mathbb{R}$. 
Then, by an induction argument on the order of differentiation  one proves that there exist real analytic functions $f_j$ from $\mathscr{E}_R(2k,n)\times\partial\mathbb{B}_n$ to $\mathbb{R}$ such that   
\begin{equation}\label{DAflatfj}
\Delta^{(n+1)/2}\bigl( |x|^{2k+1+j}A^\flat_{2k+j}(\mathbf{a},x/|x|)\bigr)= |x|^{2k-n+j}f_{j}(\mathbf{a},x/|x|)
\end{equation} for all $(\mathbf{a},x)\in \mathscr{E}_R(2k,n)\times(\mathbb{R}^n\setminus\{0\})$  and all $j\in\mathbb{N}$. Further, one has
\[
f_j(\mathbf{a},-\theta)=(-1)^j f_j(\mathbf{a},\theta)\quad\forall (\mathbf{a},\theta)\in \mathscr{E}_R(2k,n)\times\partial\mathbb{B}_n\,,\ j\in\mathbb{N}
\] and  the series
\[
\sum_{j=0}^\infty r^j f_j(\mathbf{a},\theta)
\]  converges absolutely and uniformly in the compact subsets of $\mathscr{E}_R(2k,n)\times\partial\mathbb{B}_n\times\mathbb{R}$. Hence one verifies that
\[
\begin{split}
&\Delta^{(n+1)/2}W_0(\mathbf{a},x)\\
&\quad =\Delta^{(n+1)/2}\sum_{j=0}^\infty|x|^{2k+1+j}A^\flat_{2k+j}(\mathbf{a},x/|x|)\\
&\quad =\sum_{j=0}^\infty\Delta^{(n+1)/2}\left(|x|^{2k+1+j}A^\flat_{2k+j}(\mathbf{a},x/|x|)\right)=|x|^{2k-n}\sum_{j=0}^\infty |x|^{j}f_j(\mathbf{a},x/|x|)
\end{split}
\]
for all $(\mathbf{a},x)\in \mathscr{E}_R(2k,n)\times(\mathbb{R}^n\setminus\{0\})$. Then, by equality $\Delta^{(n+1)/2}W_0(\mathbf{a},x)=|x|^{2k-n}A(\mathbf{a},x/|x|,|x|)$ and by standard properties of real analytic functions it follows that 
\[
A(\mathbf{a},\theta,r)=\sum_{j=0}^\infty r^{j}f_j(\mathbf{a},\theta)\quad \forall (\mathbf{a},\theta,r)\in \mathscr{E}_R(2k,n)\times\partial\mathbb{B}_n\times\mathbb{R}
\,.
\]
 Thus the proof for $n$ odd is complete.

Now consider the case of dimension $n$ even.  By Proposition \ref{ABC}, one has $W_1(\mathbf{a},x)=|x|^{2k}A^\sharp(\mathbf{a},x/|x|,|x|)+\log|x|\,B^\sharp(\mathbf{a},x)$ with $A^\sharp$ real analytic from $\mathscr{E}_R(2k,n)\times\partial\mathbb{B}_n\times\mathbb{R}$ to $\mathbb{R}$ and $B^\sharp$ real analytic from $\mathscr{E}_R(2k,n)\times\mathbb{R}^{n}$ to $\mathbb{R}$. If $h\in\{1,\dots,n\}$, then one verifies that 
\[
\partial_{x_h} W_{1}(\mathbf{a},x)=|x|^{2k-1}A^{\sharp,h}(\mathbf{a},x/|x|,|x|)+\log|x|\,B^{\sharp,h}(\mathbf{a},x)
\] for all $(\mathbf{a},x)\in \mathscr{E}_R(2k,n)\times(\mathbb{R}^n\setminus\{0\})$, where $A^{\sharp,h}$ is the real analytic function from $\mathscr{E}_R(2k,n)\times\partial\mathbb{B}_n\times\mathbb{R}$ to $\mathbb{R}$ defined by
\[
A^{\sharp,h}(\mathbf{a},\theta,r)\equiv 2k\theta_h A^\sharp(\mathbf{a},\theta,r)+\theta_h \,r\partial_{r}A^\sharp(\mathbf{a},\theta,r)+\mathscr{D}_{\theta_h}A^\sharp(\mathbf{a},\theta,r)+\theta_h  B^\sharp(\mathbf{a},r\theta)
\]  for all $(\mathbf{a},\theta,r)\in \mathscr{E}_R(2k,n)\times\partial\mathbb{B}_{n}\times\mathbb{R}$ and where $B^{\sharp,h}$ is the real analytic function from $\mathscr{E}_R(2k,n)\times\partial\mathbb{B}_n\times\mathbb{R}^n$ to $\mathbb{R}$ defined by
\[
B^{\sharp,h}(\mathbf{a},x)\equiv \partial_{x_h}B^{\sharp,h}(\mathbf{a},x)\qquad\forall (\mathbf{a},x)\in \mathscr{E}_R(2k,n)\times \mathbb{R}^n
\] (see also \eqref{sec2eq2} and Lemma \ref{sec2lem3}). Hence, by an induction argument on the order of differentiation  one verifies that there exist real analytic functions $A$ from $\mathscr{E}_R(2k,n)\times\partial\mathbb{B}_n\times\mathbb{R}$ to $\mathbb{R}$ and $B$ from $\mathscr{E}_R(2k,n)\times\mathbb{R}^{n}$ to $\mathbb{R}$, such that  
\[
\Delta^{n/2}W_1(\mathbf{a},x)=|x|^{2k-n}A(\mathbf{a},x/|x|,|x|)+\log|x|\,B(\mathbf{a},x)
\] for all $(\mathbf{a},x)\in \mathscr{E}_R(2k,n)\times(\mathbb{R}^n\setminus\{0\})$. Since $W_2(\mathbf{a},x)=C^\sharp(\mathbf{a},x)$  for all $(\mathbf{a},x)\in (\mathbf{a},x)\in \mathscr{E}_R(2k,n)\times(\mathbb{R}^n\setminus\{0\})$ with $C$ real analytic from $\mathscr{E}_R(2k,n)\times\mathbb{R}^{n}$ to $\mathbb{R}$ (cf.~Proposition \ref{ABC}), one deduces that equality \eqref{fsthm1eq0} for $n$ even holds with $A$, $B$ as above and $C(\mathbf{a},x)\equiv \Delta^{n/2}C^\sharp(\mathbf{a},x)$ for all $(\mathbf{a},x)\in \mathscr{E}_R(2k,n)\times \mathbb{R}^n$ (see also Theorem \ref{Sa}).

To complete the proof in the case of dimension $n$ even, one has now to show the existence of the families of functions $\{f_j\}_{j\in\mathbb{N}}$ and  $\{b_\alpha\}_{|\alpha|\ge\sup\{2k-n,0\}}$ as in the statement. By Proposition \ref{aja}, by definition \eqref{Seq2}, by the inequality in \eqref{ajless}, and by the dominated convergence theorem, one has 
\[
W_1(\mathbf{a},x)=\sum_{j=2k}^{\infty}W_{1,j}(\mathbf{a},x)\qquad\forall (\mathbf{a},x)\in\mathscr{E}_R(2k,n)\times(\mathbb{R}^n\setminus\{0\}) \,
\] with 
\[
W_{1,j}(\mathbf{a},x)\equiv -\frac{1}{(2\pi i)^{n}}\int_{\partial\mathbb{B}_n}\frac{a_j(\mathbf{a},\xi)}{j!}\ (x\cdot\xi)^{j}\ \log |x\cdot\xi|\ d\sigma_\xi
\] for all $(\mathbf{a},x)\in\mathscr{E}_R(2k,n)\times(\mathbb{R}^n\setminus\{0\})$ and $j\in\mathbb{N}$. Let $A^\sharp_j$ and $B^\sharp_j$ be the functions defined by
\begin{equation}\label{Asharpj}
A^\sharp_j(\mathbf{a},\theta)\equiv  -\frac{1}{(2\pi i)^{n}}\int_{\partial\mathbb{B}_n}\frac{a_j(\mathbf{a},\xi)}{j!}\ (\theta\cdot\xi)^{j}\ \log |\theta\cdot\xi|\ d\sigma_\xi
\end{equation} for all $(\mathbf{a},\theta)\in \mathscr{E}_R(2k,n)\times\partial\mathbb{B}_n$ and all $j\in\mathbb{N}$, and
\begin{equation}\label{Bsharpj}
B^\sharp_j(\mathbf{a},x)\equiv  -\frac{1}{(2\pi i)^{n}}\int_{\partial\mathbb{B}_n}\frac{a_j(\mathbf{a},\xi)}{j!}\ (x\cdot\xi)^{j}\ d\sigma_\xi\quad \forall(\mathbf{a},x)\in \mathscr{E}_R(2k,n)\times\mathbb{R}^n\,,\,j\in\mathbb{N}\,.
\end{equation}
Then one has
\[
W_{1,j}(\mathbf{a},x)\equiv |x|^j A^\sharp_j(\mathbf{a},x/|x|)+\log|x| B^\sharp_j(\mathbf{a},x)
\] for all $(\mathbf{a},x)\in \mathscr{E}_R(2k,n)\times(\mathbb{R}^n\setminus\{0\})$ and $j\in\mathbb{N}$. By arguing so as above for $A^\flat_j$ one proves that $A^\sharp_j$ is real analytic from $\mathscr{E}_R(2k,n)\times\partial\mathbb{B}_n$ to $\mathbb{R}$ for all $j\in\mathbb{N}$.  By Lemma \ref{fslem1} and Proposition \ref{aja} one verifies that  also $B^\sharp_j$ is real analytic from $\mathscr{E}_R(2k,n)\times\partial\mathbb{R}^n$ to $\mathbb{R}$  for all $j\in\mathbb{N}$. Further, one has
\[
A^\sharp_j(\mathbf{a},-\theta)=(-1)^{j}A^\sharp_j(\mathbf{a},\theta)\quad\forall (\mathbf{a},\theta)\in \mathscr{E}_R(2k,n)\times\partial\mathbb{B}_n\,,\, j\in\mathbb{N}\,
\] and 
\[
B^\sharp_j(\mathbf{a},tx)=t^{j}B^\sharp_j(\mathbf{a},x)\quad\forall (\mathbf{a},x)\in \mathscr{E}_R(2k,n)\times\mathbb{R}^n\,,\, t\in\mathbb{R}\,,\, j\in\mathbb{N}\,.
\] Hence one deduces that there exist real analytic functions $b^\sharp_\alpha$ from $\mathscr{E}_R(2k,n)$ to $\mathbb{R}$, for all $\alpha\in\mathbb{N}^n$, such that $B^\sharp_j(\mathbf{a},x)=\sum_{|\alpha|=j}b^\sharp_\alpha(\mathbf{a})x^\alpha$ for all $(\mathbf{a},x)\in \mathscr{E}_R(2k,n)\times\mathbb{R}^n$, $j\in\mathbb{N}$. Also note that $b^\sharp_\alpha=0$ for $|\alpha|\le 2k-1$. Then, by the inequality in \eqref{ajless} one proves that the series 
\[
\sum_{j=0}^\infty r^{j} A^\sharp_{2k+j}(\mathbf{a},\theta,r)\quad\text{and}\quad\sum_{|\alpha|\ge 2k}b^\sharp_\alpha(\mathbf{a})x^\alpha=\sum_{j=0}^\infty B^\sharp_{2k+j}(\mathbf{a},x)
\] converge absolutely and uniformly in the compact subsets of $\mathscr{E}_R(2k,n)\times\partial\mathbb{B}_n\times\mathbb{R}$ and of $\mathscr{E}_R(2k,n)\times\mathbb{R}^n$, respectively.
Now let $h\in\{1,\dots,n\}$. Then one has
\[
\begin{split}
&\partial_{x_h}\left(|x|^{2k+j} A^\sharp_{2k+j}(\mathbf{a},x/|x|,|x|)+\log|x| B^\sharp_{2k+j}(\mathbf{a},x)\right)\\
&\qquad\qquad\qquad  =|x|^{2k-1+j} A^{\sharp,h}_{2k+j}(\mathbf{a},x/|x|,|x|)+\log|x| B^{\sharp,h}_{2k+j}(\mathbf{a},x)
\end{split}
\] for all $(\mathbf{a},x)\in\mathscr{E}_R(2k,n)\times(\mathbb{R}^n\setminus\{0\})$, where $A^{\sharp,h}_j$ is the real analytic function from $\mathscr{E}_R(2k,n)\times\partial\mathbb{B}_n\times\mathbb{R}$ to $\mathbb{R}$ defined by
\[
A^{\sharp,h}_j(\mathbf{a},\theta,r)\equiv j\theta_h A^\sharp_j(\mathbf{a},\theta,r)+\theta_h \,r\partial_{r}A^\sharp_j(\mathbf{a},\theta,r)+\mathscr{D}_{\theta_h}A^\sharp_j(\mathbf{a},\theta,r)+\theta_h  B^\sharp_j(\mathbf{a},r\theta)
\]  for all $(\mathbf{a},\theta,r)\in \mathscr{E}_R(2k,n)\times\partial\mathbb{B}_{n}\times\mathbb{R}$ and all $j\in\mathbb{N}$, and where $B^{\sharp,h}_j$ is the real analytic function from $\mathscr{E}_R(2k,n)\times\partial\mathbb{B}_n\times\mathbb{R}^n$ to $\mathbb{R}$ defined by
\[
B^{\sharp,h}(\mathbf{a},x)\equiv \partial_{x_h}B^{\sharp}(\mathbf{a},x)\qquad\forall (\mathbf{a},x)\in \mathscr{E}_R(2k,n)\times\mathbb{R}^n\,,\,j\in\mathbb{N}
\] (see also \eqref{sec2eq2} and Lemma \ref{sec2lem3}). Note that 
\[
A^{\sharp,h}_j(\mathbf{a},-\theta)=(-1)^{j-1}A^{\sharp,h}_j(\mathbf{a},\theta)\quad\forall (\mathbf{a},\theta)\in \mathscr{E}_R(2k,n)\times\partial\mathbb{B}_n\,,\, j\in\mathbb{N}\,
\] and 
\[
B^{\sharp,h}_j(\mathbf{a},tx)=t^{j-1}B^{\sharp,h}_j(\mathbf{a},x)\quad\forall (\mathbf{a},x)\in \mathscr{E}_R(2k,n)\times\mathbb{R}^n\,,\, t\in\mathbb{R}\,,\, j\in\mathbb{N}\,.
\] It follows that there exist real analytic functions $b^{\sharp,h}_\alpha$ from $\mathscr{E}_R(2k,n)$ to $\mathbb{R}$, for all $\alpha\in\mathbb{N}^n$, such that $B^{\sharp,h}_{j+1}(\mathbf{a},x)=\sum_{|\alpha|=j}b^{\sharp,h}_\alpha(\mathbf{a})x^\alpha$ for all $(\mathbf{a},x)\in \mathscr{E}_R(2k,n)\times\mathbb{R}^n$, $j\in\mathbb{N}$. Also note that $b^{\sharp,h}_\alpha=0$ for $|\alpha|\le 2k-2$. Further, the series 
\[
\sum_{|\alpha|\ge 2k-1}b^{\sharp,h}_\alpha(\mathbf{a})x^\alpha=\sum_{j=0}^\infty B^{\sharp,h}_{2k+j}(\mathbf{a},x)=\partial_{x_h}\sum_{j=0}^\infty B^\sharp_{2k+j}(\mathbf{a},x)
\] converges absolutely and uniformly in the compact subsets of $\mathscr{E}_R(2k,n)\times\mathbb{R}^n$ and the series
\[
\begin{split}
&\sum_{j=0}^\infty r^{2k+j} A^{\sharp,h}_{2k+j}(\mathbf{a},\theta)\\
&\qquad=(\mathscr{D}_{\theta_h}+r\theta_h\partial_r)\sum_{j=0}^\infty r^{2k+j} A^{\sharp}_{2k+j}(\mathbf{a},\theta)+\theta_h\sum_{j=0}^\infty B^{\sharp}_{2k+j}(\mathbf{a},r\theta)
\end{split}
\] converges absolutely and uniformly in the compact subsets of $\mathscr{E}_R(2k,n)\times\partial\mathbb{B}_n\times\mathbb{R}$. So that the series 
\[
\sum_{j=0}^\infty r^{j} A^{\sharp,h}_{2k+j}(\mathbf{a},\theta,r)
\] converges absolutely and uniformly in the compact subsets of $\mathscr{E}_R(2k,n)\times\partial\mathbb{B}_n\times\mathbb{R}$. Then, by an induction argument on the order of differentiation one verifies that there exists real analytic functions $f_j$ from $\mathscr{E}_R(2k,n)\times\partial\mathbb{B}_n$ to $\mathbb{R}$ with $j\in\mathbb{N}$ and real analytic functions $b_\alpha$ from  $\mathscr{E}_R(2k,n)$ to $\mathbb{R}$ with $\alpha\in\mathbb{N}^n$, such that 
\begin {equation}\label{DABsharpj}
\begin{split}
&\Delta^{n/2}\left(|x|^{2k+j} A^\sharp_{2k+j}(\mathbf{a},x/|x|,|x|)+\log|x| B^\sharp_{2k+j}(\mathbf{a},x)\right)\\
&\qquad\qquad\qquad  =|x|^{2k-n+j} f_{j}(\mathbf{a},x/|x|)+\log|x| \sum_{|\alpha|=2k-n+j} b_{\alpha}(\mathbf{a})x^\alpha
\end{split}
\end{equation} for all $(\mathbf{a},x)\in\mathscr{E}_R(2k,n)\times(\mathbb{R}^n\setminus\{0\})$ and all $j\in\mathbb{N}$. Here we understand that $\sum_{|\alpha|=2k-n+j} b_{\alpha}(\mathbf{a})x^\alpha=0$ if $2k-n+j<0$. Further, one has
\[
f_j(\mathbf{a},-\theta)=(-1)^jf_j(\mathbf{a},\theta)\qquad\forall (\mathbf{a},\theta)\in\mathscr{E}_R(2k,n)\times\partial\mathbb{B}_n\,,\,j\in\mathbb{N}\,
\] and $b_\alpha=0$ if $|\alpha|<2k-n$. The series
\[
\sum_{j=0}^\infty r^j f_j(\mathbf{a},\theta)\quad\text{and}\quad\sum_{|\alpha|\ge\sup\{2k-n,0\}}b_\alpha(\mathbf{a})x^\alpha
\] converge absolutely and uniformly in the compact subsets of $\mathscr{E}_R(2k,n)\times\partial\mathbb{B}_n\times\mathbb{R}$ and of $\mathscr{E}_R(2k,n)\times\mathbb{R}^n$, respectively. Thus
\[
\begin{split}
&\Delta^{n/2}W_1(\mathbf{a},x)=\Delta^{n/2}\biggl(\sum_{j=0}^\infty |x|^{2k+j} A^\sharp_{2k+j}(\mathbf{a},x/|x|,|x|)+\log|x| \sum_{j=0}^\infty B^\sharp_{2k+j}(\mathbf{a},x)\biggr)\\
&\quad =\sum_{j=0}^\infty \Delta^{n/2} \left( |x|^{2k+j} A^\sharp_{2k+j}(\mathbf{a},x/|x|,|x|)+\log|x|  B^\sharp_{2k+j}(\mathbf{a},x)\right)\\
&\quad =|x|^{2k-n}\sum_{j=0}^\infty |x|^j f_j(\mathbf{a},x/|x|)+\log|x| \sum_{|\alpha|\ge\sup\{2k-n,0\}}b_\alpha(\mathbf{a})x^\alpha
\end{split}
\] for all $(\mathbf{a},x)\in\mathscr{E}_R(2k,n)\times(\mathbb{R}^n\setminus\{0\})$. Then, by equality $\Delta^{n/2}W_1(\mathbf{a},x)=|x|^{2k-n}A(\mathbf{a},x/|x|,|x|)+\log|x| B(\mathbf{a},x)$ and by standard properties of real analytic functions one verifies that 
\[
A(\mathbf{a},\theta,r)=\sum_{j=0}^\infty r^j f_j(\mathbf{a},\theta)\qquad\forall(\mathbf{a},\theta,r)\in\mathscr{E}_R(2k,n)\times\partial\mathbb{B}_n\times\mathbb{R}
\] and
\[
B(\mathbf{a},x)=\sum_{|\alpha|\ge\sup\{2k-n,0\}}b_\alpha(\mathbf{a})x^\alpha\qquad\forall(\mathbf{a},x)\in\mathscr{E}_R(2k,n)\times\mathbb{R}^n\,.
\] The proof of the theorem is now complete.\qed

\medskip

In the following Theorem \ref{fsthm2} we introduce a real analytic function $S_0$ from $\mathscr{E}_R(2k,n)\times(\mathbb{R}^n\setminus\{0\})$ to $\mathbb{R}$ such that $S_0(\mathbf{a},\cdot)$ is a fundamental solution of the principal term $\mathbf{L}_0[\mathbf{a}]$ of the operator $\mathbf{L}[\mathbf{a}]$, for all $\mathbf{a}\in\mathscr{E}_R(2k,n)$ (see also \eqref{L0}).

\begin{thm}\label{fsthm2}
Let $\{f_j\}_{j\in\mathbb{N}}$ and $\{b_\alpha\}_{|\alpha|\ge\sup\{2k-n,0\}}$  be as in Theorem \ref{fsthm1}. Let $S_{0}$ be the function from $\mathscr{E}_R(2k,n)\times(\mathbb{R}^n\setminus\{0\})$ to $\mathbb{R}$ defined by
\[
S_{0}(\mathbf{a},x)\equiv |x|^{2k-n}f_0(\mathbf{a},x/|x|)+\log |x|\sum_{|\alpha|=2k-n}b_{\alpha}(\mathbf{a})x^{\alpha}
\] for all $(\mathbf{a},x)\in\mathscr{E}_R(2k,n)\times(\mathbb{R}^n\setminus\{0\})$. Then $S_{0}(\mathbf{a},\cdot)$ is a fundamental solution of the homogeneous operator $\mathbf{L}_{0}[\mathbf{a}]$ for all fixed $\mathbf{a}\in\mathscr{E}_R(2k,n)$ (note that $\sum_{|\alpha|=2k-n}b_{\alpha}(\mathbf{a})x^{\alpha}=0$ if $n$ is odd or $\ge 2k+1$).
\end{thm}
\proof Let $v_0$ be the function from $\mathscr{E}_R(2k,n)\times\mathbb{R}^n\times\partial\mathbb{B}_n\times\mathbb{R}$ to $\mathbb{R}$ defined by
\begin{equation}\label{fsprop2eq1}
v_0(\mathbf{a},x,\xi,t)\equiv\frac{a_{2k}(\mathbf{a},\xi)}{(2k)!}(x\cdot\xi-t)^{2k}=\frac{(x\cdot\xi-t)^{2k}}{(2k)!\,P_{0}[\mathbf{a}](\xi)}
\end{equation}  for all $(\mathbf{a},x,\xi,t)\in\mathscr{E}_R(2k,n)\times\mathbb{R}^n\times\partial\mathbb{B}_n\times\mathbb{R}$ (cf.~Proposition \ref{aja}). Then one verifies that $v_0$ is the unique real analytic function from  $\mathscr{E}_R(2k,n)\times\mathbb{R}^n\times\partial\mathbb{B}_n\times\mathbb{R}$ to $\mathbb{R}$ such that $\mathbf{L}_0[\mathbf{a}]v_0(\mathbf{a},x,\xi,t)=1$ for all $(\mathbf{a},x,\xi,t)\in\mathscr{E}_R(2k,n)\times\mathbb{R}^n\times\partial\mathbb{B}_n\times\mathbb{R}$ and $\partial_x^\alpha v_0(\mathbf{a},x,\xi,t)=0$ for all $(\mathbf{a},x,\xi,t)\in \mathscr{E}_R(2k,n)\times\mathbb{R}^n\times\partial\mathbb{B}_n\times\mathbb{R}$ with $x\cdot\xi=t$ and for all $\alpha\in\mathbb{N}^n$ with $|\alpha|\le 2k-1$ (see also Proposition \ref{va}).
Let $\tilde{W}_0$, $\tilde{W}_1$, and $\tilde{W}_2$ be the functions from $\mathscr{E}_R(2k,n)\times(\mathbb{R}^n\setminus\{0\})$ to $\mathbb{R}$ defined by
\begin{eqnarray}\label{tildeW0}
&&\tilde{W}_0(\mathbf{a},x)\equiv\frac{1}{4(2\pi i)^{n-p_n}}\int_{\partial\mathbb{B}_n}\int_0^{x\cdot\xi}v_0(\mathbf{a},x,\xi,t)\,\mathrm{sgn}\,t\ dt\,d\sigma_\xi\,,\\
\label{tildeW1}
&&\tilde{W}_1(\mathbf{a},x)\equiv-\frac{1}{(2\pi i)^{n-p_n}}\int_{\partial\mathbb{B}_n} v_0(\mathbf{a},x,\xi,0) \log|x\cdot \xi|\, d\sigma_\xi\,,\\
\label{tildeW2}
&&\tilde{W}_2(\mathbf{a},x)\equiv-\frac{1}{(2\pi i)^{n-p_n}}\int_{\partial\mathbb{B}_n} \int_0^{x\cdot\xi}\frac{v_0(\mathbf{a},x,\xi,t)-v_0(\mathbf{a},x,\xi,0)}{t}\ dt\,d\sigma_\xi\qquad
\end{eqnarray} for all $(\mathbf{a},x)\in\mathscr{E}_R(2k,n)\times(\mathbb{R}^n\setminus\{0\})$. If the dimension $n$ is odd,  then Theorem \ref{Sa} implies that the function  $\Delta^{(n+1)/2}\tilde{W}_0(\mathbf{a},\cdot)$ is a fundamental solution of $\mathbf{L}_{0}[\mathbf{a}]$, for all $\mathbf{a}\in\mathscr{E}_R(2k,n)$. Moreover, by  
\eqref{Aflatj}, \eqref{DAflatfj}, \eqref{fsprop2eq1}, \eqref{tildeW0}, and by a straightforward calculation, one verifies that $\Delta^{(n+1)/2}\tilde{W}_0(\mathbf{a},x)=|x|^{2k-n}f_0(\mathbf{a},x/|x|)$ for all $(\mathbf{a},x)\in\mathscr{E}_R(2k,n)\times(\mathbb{R}^n\setminus\{0\})$.   Hence the validity of the theorem for $n$ odd is proved.  Now let $n$ be even.  Then one can verify that $\tilde{W}_2(\mathbf{a},\cdot)$ equals a polynomial function of degree $2k$ on $\mathbb{R}^n$, for all fixed $\mathbf{a}\in\mathscr{E}_R(2k,n)$. Thus
$\mathbf{L}_{0}[\mathbf{a}]\Delta^{n/2}\tilde{W}_2(\mathbf{a},\cdot)=0$ because the operator $\mathbf{L}_{0}[\mathbf{a}]\Delta^{n/2}$ is homogeneous of order $2k+n> 2k+1$.
Hence, Theorem \ref{Sa}   implies that the function   $\Delta^{n/2}\tilde{W}_1(\mathbf{a},\cdot)$  is a fundamental solution of $\mathbf{L}_0[\mathbf{a}]$ for all $\mathbf{a}\in\mathscr{E}_R(2k,n)$. Moreover, by  \eqref{Asharpj},  \eqref{Bsharpj}, \eqref{DABsharpj}, \eqref{fsprop2eq1}, \eqref{tildeW1}, and by a straightforward calculation, one verifies that $\Delta^{n/2}\tilde{W}_1(\mathbf{a},x)=|x|^{2k-n}f_0(\mathbf{a},x/|x|)+\log |x|\sum_{|\alpha|=2k-n}b_{\alpha}(\mathbf{a})x^{\alpha}$ for all $(\mathbf{a},x)\in\mathscr{E}_R(2k,n)\times(\mathbb{R}^n\setminus\{0\})$.  Now the proof is complete. 
\qed\medskip

\section{The  single layer potential $v[\mathbf{a},\mu]$}\label{slsec}

In this section we show some properties of the single layer  potential $v[\mathbf{a},\mu]$ corresponding to the fundamental solution $S(\mathbf{a},\cdot)$ introduced in Theorem \ref{fsthm1} (see also definition \eqref{sl}). 
We introduce the following notation. Let $\mathbf{a}\in\mathscr{E}_R(2k,n)$. Let $\mu\in C^{m-1,\lambda}(\partial\Omega)$.
Let $\beta\in\mathbb{N}^n$ and $|\beta|\le 2k-1$. Then $v_\beta[\mathbf{a},\mu]$ denotes the function from $\mathbb{R}^n$ to $\mathbb{R}$ defined by  
\[
v_\beta[\mathbf{a},\mu](x)\equiv\int_{\partial\Omega}\partial_x^\beta S(\mathbf{a},x-y)\mu(y)\,d\sigma_y\qquad\forall x\in\mathbb{R}^n\,,
\] where the integral is understood in the sense of singular integrals if $x\in\partial\Omega$ and $|\beta|=2k-1$. Thus, for $\beta=(0,\dots,0)$ one has 
\[
v_{(0,\dots,0)}[\mathbf{a},\mu]= v[\mathbf{a},\mu]\,.
\]  Moreover, by standard theorems on differentiation under integral sign one verifies that 
\begin{equation}\label{2k-1}
v_\beta[\mathbf{a},\mu](x)=\partial_x^\beta v[\mathbf{a},\mu](x)\qquad\forall x\in\mathbb{R}^n\setminus\partial\Omega\,,\beta\in\mathbb{N}^n\,,\, |\beta|\le 2k-1\,.
\end{equation}
Then, the following Theorem \ref{2k-2}  implies that $v[\mathbf{a},\mu]\in C^{2k-2}(\mathbb{R}^n)$ and that $\partial^\beta_x v[\mathbf{a},\mu]=v_\beta[\mathbf{a},\mu]$ in the whole of $\mathbb{R}^n$ for all $\beta\in\mathbb{N}^n$ with $|\beta|\le 2k-2$.

\begin{thm}\label{2k-2}
 Let $\mathbf{a}\in\mathscr{E}_R(2k,n)$. Let $\mu\in C^{m-1,\lambda}(\partial\Omega)$. Let $\beta\in\mathbb{N}^n$, $|\beta|\le 2k-2$. Then, $v_\beta[\mathbf{a},\mu]\in C^{2k-2-|\beta|}(\mathbb{R}^n)$ and $\partial_x^\beta v[\mathbf{a},\mu](x)=v_{\beta}[\mathbf{a},\mu](x)$ for all $x\in\mathbb{R}^n$. \end{thm}
\proof Let $B$ be the function in Theorem \ref{fsthm1}. By standard properties of real analytic functions one verifies that the function from $\mathscr{E}_R(2k,n)\times\partial\mathbb{B}_n\times(\mathbb{R}\setminus\{0\})$ to $\mathbb{R}$ which takes $(\mathbf{a},\theta,r)$ to $r^{-(2k-n)}B(\mathbf{a},r\theta)$ extends to a unique real analytic function $\tilde{B}$ from $\mathscr{E}_R(2k,n)\times\partial\mathbb{B}_n\times\mathbb{R}$ to $\mathbb{R}$. Then by Theorem \ref{fsthm1} one has
\[
S(\mathbf{a},x)=|x|^{2k-n}A(\mathbf{a},x/|x|,|x|)+|x|^{2k-n}\log |x|\tilde{B}(\mathbf{a},x/|x|,|x|)+C(\mathbf{a},x)
\] for all $(\mathbf{a},x)\in\mathscr{E}_R(2k,n)\times(\mathbb{R}^n\setminus\{0\})$. Now let $h\in\{1,\dots,n\}$.  Then a straightforward calculation shows that 
\[
\begin{split}
&\partial_{x_h} S(\mathbf{a},x)\\
&\ =|x|^{2k-n-1}A_h(\mathbf{a},x/|x|,|x|)+|x|^{2k-n-1}\log |x|\tilde{B}_h(\mathbf{a},x/|x|,|x|)+\partial_{x_h} C(\mathbf{a},x)
\end{split}
\] for all $(\mathbf{a},x)\in\mathscr{E}_R(2k,n)\times(\mathbb{R}^n\setminus\{0\})$, where $A_h$ and  $\tilde{B}_h$ are the real analytic functions from $\mathscr{E}_R(2k,n)\times\partial\mathbb{B}_n\times\mathbb{R}$ to $\mathbb{R}$ defined by 
\[
\begin{split}
&A_h(\mathbf{a},\theta,r)\equiv\mathscr{D}_{\theta_h} A(\mathbf{a},\theta,r)+r\theta_h\partial_r A(\mathbf{a},\theta,r)+(2k-n)\theta_h A(\mathbf{a},\theta,r)+\theta_h  \tilde{B}(\mathbf{a},\theta,r)\,,\\
&\tilde{B}_h(\mathbf{a},\theta,r)\equiv\mathscr{D}_{\theta_h}\tilde{B}(\mathbf{a},\theta,r)+r\theta_h\partial_r\tilde{B}(\mathbf{a},\theta,r)+(2k-n)\theta_h \tilde{B}(\mathbf{a},\theta,r)\,,
\end{split}
\] for all $(\mathbf{a},\theta,r)\in\mathscr{E}_R(2k,n)\times\partial\mathbb{B}_n\times\mathbb{R}$ (see also Lemma \ref{sec2lem3} and equality \eqref{sec2eq2}). Then, by an induction argument on the order of differentiation one verifies that for all $\alpha\in\mathbb{N}^n$ there exist real analytic functions $A_\alpha$ and  $\tilde{B}_\alpha$  from $\mathscr{E}_R(2k,n)\times\partial\mathbb{B}_n\times\mathbb{R}$ to $\mathbb{R}$ such that  
\[
\begin{split}
&\partial_x^\alpha S(\mathbf{a},x)\\
&\ =|x|^{2k-n-|\alpha|}A_\alpha(\mathbf{a},x/|x|,|x|)+|x|^{2k-n-|\alpha|}\log |x|\tilde{B}_\alpha(\mathbf{a},x/|x|,|x|)+\partial_x^\alpha C(\mathbf{a},x)
\end{split}
\] for all $(\mathbf{a},x)\in\mathscr{E}_R(2k,n)\times(\mathbb{R}^n\setminus\{0\})$. In particular, $\partial^\alpha_x S(\mathbf{a},x)=o(|x|^{2-n-1/2})$ as $|x|\to 0^+$ for all $\alpha\in\mathbb{N}^n$ with $|\alpha|\le 2k-2$. Then, by the Vitali Convergence Theorem one can prove that the function $v_{\alpha}[\mathbf{a},\mu]$ is continuous on $\mathbb{R}^n$ for all $\alpha\in\mathbb{N}^n$, $|\alpha|\le 2k-2$  (see also  Folland \cite[Proposition 3.25]{F} where the continuity of the single layer potential corresponding to the fundamental solution of the Laplace operator is proved, but the proof for $v_{\alpha}[\mathbf{a},\mu]$ is based on the same argument). Then, by the Divergence Theorem one verifies that
\[
\begin{split}
&\int_{\mathbb{R}^n}\phi(x)v_{\alpha}[\mathbf{a},\mu](x)\,dx\\
&\qquad=\int_{\Omega}\phi(x)v_{\alpha}[\mathbf{a},\mu](x)\,dx+\int_{\mathbb{R}^n\setminus\mathrm{cl}\Omega}\phi(x)v_{\alpha}[\mathbf{a},\mu](x)\,dx\\
&\qquad=\int_{\mathbb{R}^n}(\partial_x^\alpha\phi(x)) v[\mathbf{a},\mu](x)\,dx
\end{split}
\] for all smooth function $\phi$ from $\mathbb{R}^n$ to $\mathbb{R}$ with compact support (see also \eqref{2k-1}). Hence, $\partial_x^{\alpha}v[\mathbf{a},\mu]=v_{\alpha}[\mathbf{a},\mu]$ in the sense of distributions on $\mathbb{R}^n$. Then a standard argument based on the convolution with a family of mollifiers shows that $v[\mathbf{a},\mu]\in C^{2k-2}(\mathbb{R}^n)$ and that  $\partial_x^{\alpha}v[\mathbf{a},\mu](x)=v_{\alpha}[\mathbf{a},\mu](x)$ for all $x\in\mathbb{R}^n$  and for all $\alpha\in\mathbb{N}^n$ with $|\alpha|\le 2k-2$. Thus $\partial_x^\beta v[\mathbf{a},\mu]\in C^{2k-2-|\beta|}(\mathbb{R}^n)$ and the validity of the theorem is proved.\qed\medskip

As we have seen in Theorem \ref{2k-2}, $v[\mathbf{a},\mu]$ is a function of class $C^{2k-2}$.  In Theorem \ref{fs-vs-m} here below we show that the restriction of $v[\mathbf{a},\mu]$ to $\Omega$ can be extended to a function of class $C^{m+2k-2,\lambda}$ from $\mathrm{cl}\Omega$ to $\mathbb{R}$ (a similar result  holds for the restriction of $v[\mathbf{a},\mu]$ to the exterior domain $\mathbb{R}^n\setminus\mathrm{cl}\Omega$ in a local sense which will be clarified). In order to prove Theorem \ref{fs-vs-m} we exploit an idea of Miranda (cf.~\cite[\S 5]{Mi65}). Accordingly, we first state in Theorem \ref{2.1} a result by Miranda (cf.~\cite[Theorem 2.1]{Mi65}). 

\begin{thm}
\label{2.1}  Let $K\in C^{2m}(\mathbb{R}^n\setminus\{0\})$ be such that 
\begin{equation}\label{2.1eq1} 
K(x)=|x|^{1-n}K(x/|x|)\qquad \forall x\in\mathbb{R}^n\setminus\{0\}\,,
\end{equation} and 
\begin{equation}\label{2.1eq2}
\int_{\Pi\cap\partial\mathbb{B}_n} K(\eta)\ d\sigma_\eta=0\text{ for every hyper-plane $\Pi$ of $\mathbb{R}^n$ which contains $0$.}
\end{equation}      Let $p[\mu]$ be the function from $\mathbb{R}^n\setminus\partial\Omega$ to $\mathbb{R}$ defined by
\[
p[\mu](x)\equiv\int_{\partial\Omega}K(x-y)\mu(y)\, d\sigma_y\qquad\forall x\in\mathbb{R}^n\setminus\partial\Omega
\] for all $\mu\in C^{m-1,\lambda}(\partial\Omega)$.
Then the following statements hold.
\begin{enumerate}
\item[(i)] The restriction $p[\mu]_{|\Omega}$ extends to a unique continuous function $p^+[\mu]$ on $\mathrm{cl}\Omega$ for all $\mu\in C^{m-1,\lambda}(\partial\Omega)$. The map which takes $\mu$ to $p^+[\mu]$ is continuous from $C^{m-1,\lambda}(\partial\Omega)$ to $C^{m-1,\lambda}(\mathrm{cl}\Omega)$.
\item[(ii)] The restriction $p[\mu]_{|\mathbb{R}^n\setminus\mathrm{cl}\Omega}$ extends to a unique continuous function $p^-[\mu]$ on $\mathbb{R}^n\setminus\Omega$ for all $\mu\in C^{m-1,\lambda}(\partial\Omega)$. If $R\in]0,+\infty[$ and $\mathrm{cl}\Omega$ is contained in $R\mathbb{B}_{n}$, then the map which takes $\mu$ to $p^-[\mu]_{|\mathrm{cl}(R\mathbb{B}_{n})\setminus\Omega}$ is continuous from $C^{m-1,\lambda}(\partial\Omega)$ to $C^{m-1,\lambda}(\mathrm{cl}(R\mathbb{B}_{n})\setminus\Omega)$.
\end{enumerate}
\end{thm}

We now introduce in the following Lemmas \ref{fs-vs-m lemma1}--\ref{fs-vs-m lemma5} some technical results which are needed in the proof of Theorem \ref{fs-vs-m}.

\begin{lem}
\label{fs-vs-m lemma1}
Let $j,q\in\mathbb{N}$, $q\geq 1$.  Let $\alpha\in\mathbb{N}^n$ and $|\alpha|=q-1$. Let $f$ be a real analytic function from $\partial\mathbb{B}_n$ to $\mathbb{R}$  such that $f(-\theta)=(-1)^jf(\theta)$ for all $\theta\in\partial\mathbb{B}_n$. 
 Let $K(x)\equiv\partial_x^{\alpha}(|x|^{q-n}f(x/|x|))$ for all $x\in\mathbb{R}^n\setminus\{0\}$. If the sum $j+q$ is even, then  $K$ satisfies the conditions in \eqref{2.1eq1} and \eqref{2.1eq2}.
\end{lem}
\proof The lemma clearly holds for $q=1$ and $j$ odd. Let $q\ge 2$.  Then, one has $\partial_{x_h}(|x|^{q-n} f(x/|x|))=|x|^{q-n-1}g_h(x/|x|)$ for all $x\in\mathbb{R}^n\setminus\{0\}$ and all $h\in\{1,\dots,n\}$, where $g_h$ is the real analytic function from $\partial\mathbb{B}_n$ to $\mathbb{R}$ defined by $g_h(\theta)\equiv \mathscr D_{\theta_h} f(\theta)+(q-n)\theta_h f(\theta)$ for all $\theta\in \partial\mathbb{B}_n$ and $h\in\{1,\dots,n\}$ (see also Lemma \ref{sec2lem3} and equality \eqref{sec2eq2}). Moreover,  one has  $g_h(-\theta)=(-1)^{j+1}g_h(\theta)$ for all $\theta\in \partial\mathbb{B}_n$. Then,  by an induction argument on the order of differentiation  one shows that there exist real analytic functions $g_\beta$ from $\partial\mathbb{B}_n$ to $\mathbb{R}$ such that $\partial_x^{\beta}(|x|^{q-n} f(x/|x|))=|x|^{q-n-|\beta|}g_\beta(x/|x|)$ for all $x\in\mathbb{R}^n\setminus\{0\}$ and for all $\beta\in\mathbb{N}^n$. Moreover, $g_{\beta}(-\theta)=(-1)^{j+|\beta|}g_\beta(\theta)$ for all $\theta\in \partial\mathbb{B}_n$ and for all $\beta\in\mathbb{N}^n$. Now the validity of the lemma follows by taking $\beta=\alpha$, by assumption $|\alpha|=q-1$, and by a straightforward verification.\qed

\begin{lem}\label{fs-vs-m lem3}
Let $j\in\mathbb{N}$. Let $\alpha\in\mathbb{N}^n$ and $|\alpha|=n+j-1$. Let $p$ be a real homogeneous polynomial function from $\mathbb{R}^n$ to $\mathbb{R}$ of degree $j$.  
Let $K(x)\equiv\partial_x^{\alpha}(p(x)\,\log |x|)$ for all $x\in\mathbb{R}^n\setminus\{0\}$. If the dimension $n$ is even, then $K$ satisfies the conditions in \eqref{2.1eq1} and \eqref{2.1eq2}.
\end{lem}
\proof One has
\begin{equation}\label{fs-vs-m lem3eq1}
\begin{split}
&\partial_x^{\alpha}\bigl(p(x)\,\log |x|\bigr)=\sum_{\alpha'\leq\alpha}\binom{\alpha}{\alpha'}\left(\partial_x^{\alpha-\alpha'}p(x)\right)\left(\partial_x^{\alpha'}\log |x|\right)\\
&=\left(\partial_x^{\alpha}p(x)\right)\,\log |x|\,+\sum_{0<\alpha'\leq\alpha}\binom{\alpha}{\alpha'}\left(\partial_x^{\alpha-\alpha'}p(x)\right)\left(\partial_x^{\alpha'}\log |x|\right)\quad \forall x\in\mathbb{R}^n\setminus\{0\}\,.
\end{split}
\end{equation} Let now $\beta\in\mathbb{N}^n$. If $|\beta|\le j$, then $\partial_x^{\beta}p$ is a homogeneous polynomial function of degree $j-|\beta|$.  If instead $|\beta|> j$, then $\partial_x^{\beta}p=0$. Moreover, if $|\beta|>0$, then there exist $\beta',\iota\in\mathbb{N}^n$ with $|\iota|=1$ such that $\beta=\beta'+\iota$. Thus $\partial_x^{\beta}\log |x|=\partial_x^{\beta'}(\partial_x^{\iota}\log |x|)=\partial_x^{\beta'}(|x|^{-1}(x/|x|)^\iota)$. Then, by arguing so as in the proof of  Lemma \ref{fs-vs-m lemma1} it follows that  there exists a real analytic function $h_\beta$ from $\partial\mathbb{B}_n$ to $\mathbb{R}$ such that such that $\partial_x^{\beta}\log|x|=|x|^{-|\beta|} h_\beta(x/|x|)$ for all $x\in\mathbb{R}^n\setminus\{0\}$ and such that $h_\beta(-\theta)=(-1)^{|\beta|}h_\beta(\theta)$ for all $\theta\in\partial\mathbb{B}_n$. Thus equality \eqref{fs-vs-m lem3eq1} implies that 
$\partial_x^{\alpha}(p(x)\,\log |x|)={q}(x)\log |x|\,+|x|^{j-|\alpha|}g(x/|x|)$  for all $x\in\mathbb{R}^n\setminus\{0\}$ where $q$  is a real homogeneous polynomial function of degree $j-|\alpha|$ from $\mathbb{R}^n$ to $\mathbb{R}$ and $g$ is real analytic function from $\partial\mathbb{B}_{n}$ to $\mathbb{R}$ such that $g(-\theta)=(-1)^{j-|\alpha|}g(\theta)$ for all $\theta\in\partial\mathbb{B}_{n}$. Now, by the assumption $|\alpha|=n+j-1$ and by a straightforward verification the validity of the lemma follows.\qed

\medskip

\begin{lem}\label{fs-vs-m lemma5}
Let $j\in\mathbb{N}$, $j\ge 1$. Let $F$, $G$ be real analytic functions from $\partial\mathbb{B}_n\times\mathbb{R}$ to $\mathbb{R}$. Let $H$ be the function from $\mathbb{R}^n\setminus\{0\}$ to $\mathbb{R}$ defined by $H(x)\equiv |x|^{j} F(x/|x|,|x|) + |x|^{j}\log |x| G(x/|x|,|x|)$ for all $x\in\mathbb{R}^n\setminus\{0\}$. Then $H$ extends to an element of $C^{j-1}(\mathbb{R}^n)$.
\end{lem} 
\proof The statement clearly holds for $j=1$. Let $j\ge 2$. Then one notes that 
\[
\partial_{x_h}H(x)=|x|^{j-1}{F}_h(x/|x|,|x|)+|x|^{j-1}\log |x| {G}_h(x/|x|,|x|)\qquad\forall x\in\mathbb{R}^n\setminus\{0\}
\] for all $l\in\{1,\dots,n\}$, where ${F}_h$ and ${G}_h$ are real analytic functions from $\partial\mathbb{B}_n\times\mathbb{R}$ to $\mathbb{R}$. Then by arguing by induction on $j$ one proves the validity of the lemma.\qed

\medskip

We are now ready to state and  prove Theorem \ref{fs-vs-m}.

\begin{thm}\label{fs-vs-m}
Let $\mathbf{a}\in\mathscr{E}_R(2k,n)$.  Then the following statements hold.
\begin{itemize} 
\item[$(i)$] The map    which takes $\mu$ to $v[\mathbf{a},\mu]_{|\mathrm{cl}\Omega}$ is  continuous from  $C^{m-1,\lambda}(\partial\Omega)$ to $C^{m+2k-2,\lambda}(\mathrm{cl}\Omega)$. 
\item[$(ii)$] If $R>0$ is such that $\mathrm{cl}\Omega\subseteq R\mathbb{B}_{n}$, then the map which takes $\mu$ to $v[\mathbf{a},\mu]_{|\mathrm{cl}(R\mathbb{B}_n)\setminus\Omega}$ is continuous from $C^{m-1,\lambda}(\partial\Omega)$ to $C^{m+2k-2,\lambda}(\mathrm{cl}(R\mathbb{B}_n)\setminus\Omega)$.
\end{itemize}
\end{thm}
\proof Let $A$, $B$, $C$, $\{f_j\}_{j\in\mathbb{N}}$, $\{b_\alpha\}_{|\alpha|\ge\sup\{2k-n,0\}}$ be as in Theorem \ref{fsthm1}. Then
\begin{equation}\label{fs-vs-m1}
\begin{split}
&S(\mathbf{a},x)=\sum_{j=0}^{m+n-1}f_{j}(\mathbf{a},x/|x|)|x|^{2k-n+j}+\sum_{2k-n\le|\alpha|\le m+2k-1}b_\alpha(\mathbf{a}) x^\alpha\log |x|\\ 
&\quad +|x|^{m+2k} {A}_{\mathbf{a},m}(x/|x|,|x|)+|x|^{m+2k}\log |x|{B}_{\mathbf{a},m}(x/|x|,|x|)+C(\mathbf{a},x)
\end{split}
\end{equation}  for all $x\in\mathbb{R}^n\setminus\{0\}$, where ${A}_{\mathbf{a},m}$ and ${B}_{\mathbf{a},m}$ are real analytic functions from $\partial\mathbb{B}_n\times\mathbb{R}$ to $\mathbb{R}$. \par
 Note that the terms $f_j(\mathbf{a},x/|x|)|x|^{2k-n+j}$ which appear in \eqref{fs-vs-m1} are functions of the form considered in Lemma \ref{fs-vs-m lemma1}. Accordingly, the function which takes $x\in\mathbb{R}^n\setminus\{0\}$ to $\partial_x^\beta(f_j(\mathbf{a},x/|x|)|x|^{2k-n+j})$ satisfies the conditions in \eqref{2.1eq1} and \eqref{2.1eq2} for all $\beta\in\mathbb{N}^n$ with $|\beta|=2k-1+j$. Hence Theorem \ref{2.1} implies that the map from $C^{m-1,\lambda}(\partial\Omega)$ to $C^{m-1,\lambda}(\mathrm{cl}\Omega)$ which takes $\mu$ to  the unique extension to $\mathrm{cl}\Omega$ of 
\[
\int_{\partial\Omega}\partial_x^\beta\left(f_j\left(\mathbf{a},(x-y)/|x-y|\right)|x-y|^{2k-n+j}\right)\mu(y)\ d\sigma_y\qquad\forall x\in\Omega
\]  is continuous for all $\beta\in\mathbb{N}^n$ with $|\beta|=2k-1+j$. Then by the continuous embedding of $C^{m+2k-2+j,\lambda}(\mathrm{cl}\Omega)$ in $C^{m+2k-2,\lambda}(\mathrm{cl}\Omega)$ and by the equality 
\[
\begin{split}
&\partial_x^\beta\int_{\partial\Omega}f_j\left(\mathbf{a},(x-y)/|x-y|\right)|x-y|^{2k-n+j}\mu(y)\ d\sigma_y\\
&\quad=\int_{\partial\Omega}\partial_x^\beta\left(f_j\left(\mathbf{a},(x-y)/|x-y|\right)|x-y|^{2k-n+j}\right)\mu(y)\ d\sigma_y\qquad\forall x\in\Omega
\end{split}
\] one deduces that the map from $C^{m-1,\lambda}(\partial\Omega)$ to $C^{m+2k-2,\lambda}(\mathrm{cl}\Omega)$ which takes $\mu$ to  the unique extension to $\mathrm{cl}\Omega$ of 
\[
\int_{\partial\Omega} f_j\left(\mathbf{a},(x-y)/|x-y|\right)|x-y|^{2k-n+j}\mu(y)\ d\sigma_y\qquad\forall x\in\Omega
\]  is continuous. Similar result one has if $\mathrm{cl}\Omega$ is replaced by $\mathrm{cl}(R\mathbb{B}_n)\setminus\Omega$.\par 
Now consider the terms $b_{\alpha}(\mathbf{a})x^\alpha\log |x|$ which appear in \eqref{fs-vs-m1}. By Lemma \ref{fs-vs-m lem3} one verifies that $\partial_x^\beta(b_{\alpha}(\mathbf{a})x^\alpha\log |x|)$ satisfies the conditions in \eqref{2.1eq1} and \eqref{2.1eq2} for all $\beta\in\mathbb{N}^n$ with $|\beta|=n+|\alpha|-1$. Then,  by the continuous embedding of $C^{m-1+(n+|\alpha|-1),\lambda}(\mathrm{cl}\Omega)$ in $C^{m+2k-2,\lambda}(\mathrm{cl}\Omega)$ and by arguing so as above for the terms $f_j(\mathbf{a},x/|x|)|x|^{2k-n+j}$, one can prove that the map from $C^{m-1,\lambda}(\partial\Omega)$ to  $C^{m+2k-2,\lambda}(\mathrm{cl}\Omega)$ which takes $\mu$ to the the unique extension to $\mathrm{cl}\Omega$ of 
\[
\int_{\partial\Omega}b_\alpha(\mathbf{a})\,(x-y)^\alpha\log|x-y|\,\mu(y) d\sigma_y\qquad\forall x\in\Omega
\] is continuous.   Similar result one has if $\mathrm{cl}\Omega$ is replaced by $\mathrm{cl}{(R\mathbb{B}_n)\setminus\Omega}$.\par 
Finally note that the term 
\[
|x|^{m+2k}{A}_{\mathbf{a},m}(x/|x|,|x|)+|x|^{m+2k}\log |x|{B}_{\mathbf{a},m}(x/|x|,|x|)
\] in \eqref{fs-vs-m1} is a function of the form considered in Lemma \ref{fs-vs-m lemma5}  and hence extends to an element of $C^{m+2k-1}(\mathbb{R}^n)$. Moreover $C(\mathbf{a},\cdot)$ is real analytic from $\mathbb{R}^n$ to $\mathbb{R}$ and thus in $C^{m+2k-1}(\mathbb{R}^n)$.
\par
Thus, by equality \eqref{fs-vs-m1}, and by standard theorems on differentiation under integral sign, and by  the continuity of the embedding of $C^{m+2k-1}(\mathbb{R}^n)$ into $C^{m+2k-2,\lambda}(\mathbb{R}^n)$, one completes the proof of the theorem. 
\qed\medskip

As an immediate consequence of Theorem \ref{fs-vs-m} and of equality \eqref{2k-1}, one verifies the validity of the following Corollary \ref{fs-vs-m-cor}. 

\begin{cor}\label{fs-vs-m-cor}
 Let $\mathbf{a}\in\mathscr{E}_R(2k,n)$. Let $\mu\in C^{m-1,\lambda}(\partial\Omega)$. Let $\beta\in\mathbb{N}^n$, $|\beta|=2k-1$. Then the following statements hold.
\begin{itemize}
\item [$(i)$] The restriction $v_\beta[\mathbf{a},\mu]_{|\Omega}$ extends to a unique continuous function $v_\beta^+[\mathbf{a},\mu]$ on $\mathrm{cl}\Omega$ which belongs to $C^{m-1,\lambda}(\mathrm{cl}\Omega)$. 
\item[$(ii)$] The restriction $v_\beta[\mathbf{a},\mu]_{|\mathbb{R}^n\setminus\mathrm{cl}\Omega}$ extends to a unique continuous function $v_\beta^-[\mathbf{a},\mu]$ on $\mathbb{R}^n\setminus\Omega$ and $v_\beta^-[\mathbf{a},\mu]_{|\mathrm{cl}(R\mathbb{B}_n)\setminus\Omega}\in C^{m-1,\lambda}(\mathrm{cl}(R\mathbb{B}_n)\setminus\Omega)$ for all $R>0$ such that $\mathrm{cl}\Omega\subseteq R\mathbb{B}_{n}$.
\end{itemize}
\end{cor}

Finally, we consider in the following Theorem \ref{jump} the jump properties of the single layer potential $v[\mathbf{a},\mu]$.

\begin{thm}\label{jump}
Let $\mathbf{a}\in\mathscr{E}(2k,n)$. Let $\mu\in C^{m-1,\lambda}(\partial\Omega)$. Let $\beta\in\mathbb{N}^n$ and $|\beta|=2k-1$. Let $v_\beta^+[\mathbf{a},\mu]$ and $v_\beta^-[\mathbf{a},\mu]$ be as in Corollary \ref{fs-vs-m-cor}.  Then
\begin{equation}\label{jumpeqn1}
v_{\beta}^{\pm}[\mathbf{a},\mu](x) =\mp\frac{\nu_{\Omega}(x)^\beta}{2P_{0}[\mathbf{a}](\nu_{\Omega}(x))}\
 \mu(x)+v_{\beta}[\mathbf{a},\mu](x)\qquad\forall x\in\partial\Omega\,,
\end{equation} where $\nu_\Omega$ denotes the outward unit normal to the boundary of $\Omega$.
\end{thm}
\proof Let ${S}_{0}(\mathbf{a},\cdot)$ be the function in Theorem~\ref{fsthm1} and 
\[
{v}_{{0},\beta}[\mathbf{a},\mu](x)\equiv \int_{\partial\Omega}\partial^{\beta}_{x}S_{0}(\mathbf{a},x-y)\mu(y)\, d\sigma_y\qquad\forall x\in\mathbb{R}^n\setminus\partial\Omega\,.
\] Since ${S}_{0}(\mathbf{a},\cdot)$ is a fundamental solution of the homogeneous operator $\mathbf{L}_0[\mathbf{a}]$, Corollary \ref{fs-vs-m-cor} implies that  $v_{0,\beta}[\mathbf{a},\mu]_{|\Omega}$ and $v_{0,\beta}[\mathbf{a},\mu]_{|\mathbb{R}^n\setminus\mathrm{cl}\Omega}$ extend to unique continuous functions $v_{0,\beta}^+[\mathbf{a},\mu]$ on $\mathrm{cl}\Omega$ and  $v_{0,\beta}^-[\mathbf{a},\mu]$ on $\mathbb{R}^n\setminus\Omega$, respectively. Moreover,
\begin{equation}\label{jumpeqn2}
v_{0,\beta}^{\pm}[\mathbf{a},\mu](x) =\mp\frac{\nu_{\Omega}(x)^\beta}{2P_{0}[\mathbf{a}](\nu_{\Omega}(x))}\
 \mu(x)+v_{0,\beta}[\mathbf{a},\mu](x)\qquad\forall x\in\partial\Omega
\end{equation} (cf., {\it e.g.}, Mitrea \cite[pp.~392--393]{Mit10}, see also Cialdea \cite[\S2, IX]{Cia95}, \cite[Theorem~3]{Cia07}). Now let 
\[
S_\infty(\mathbf{a},x)\equiv S(\mathbf{a},x)-S_{0}(\mathbf{a},x)\qquad\forall x\in\mathbb{R}^n\setminus\{0\}
\] and define 
\[
v_{\infty,\beta}[\mathbf{a},\mu](x)\equiv \int_{\partial\Omega}\partial^{\beta}_{x}S_\infty(\mathbf{a},x-y)\mu(y)\, d\sigma_y\qquad\forall x\in\mathbb{R}^n\setminus\partial\Omega\,.
\] Thus 
\begin{equation}\label{jumpeqn3}
{v}_{\beta}[\mathbf{a},\mu]={v}_{0,\beta}[\mathbf{a},\mu]+ v_{\infty,\beta}[\mathbf{a},\mu]\,.
\end{equation}  By Theorem \ref{fsthm1} and by a straightforward verification, one has
\[
\begin{split}
&S_\infty(\mathbf{a},x)\\
&\quad =|x|^{2k+1-n}A_{\infty}(\mathbf{a},x/|x|,|x|)+|x|^{2k+1-n}\log |x| B_{\infty}(x/|x|,|x|)+C(\mathbf{a},x)
\end{split}
\] for all $x\in\mathbb{R}^n\setminus\{0\}$, where $A_{\infty}$ and $B_{\infty}$  are real analytic functions from  $\mathscr{E}_R(2k,n)\times\partial\mathbb{B}_n\times\mathbb{R}$ to $\mathbb{R}$.
Then, by arguing so as in the proof  Theorem \ref{2k-2} one proves that 
\[
\begin{split}
&\partial_{x}^\beta S_\infty(\mathbf{a},x)\\
&\quad =|x|^{2-n}A_{\infty,\beta}(\mathbf{a},x/|x|,|x|)+|x|^{2-n}\log |x|A_{\infty,\beta}(\mathbf{a},x/|x|,|x|)
+\partial_{x}^\beta C(\mathbf{a},x)
\end{split}
\]
for all $x\in\mathbb{R}^n\setminus\{0\}$, where $A_{\infty,\beta}$ and $B_{\infty,\beta}$ are  real analytic functions from  $\mathscr{E}_R(2k,n)\times\partial\mathbb{B}_n\times\mathbb{R}$ to $\mathbb{R}$. Then, by the Vitali Convergence Theorem one proves that the function $v_{\infty,\beta}[\mathbf{a},\mu]$ is continuous on $\mathbb{R}^n$   (see also  Folland \cite[Proposition 3.25]{F} where the continuity of the single layer potential corresponding to the fundamental solution of the Laplace operator is proved, but the proof for $v_{\infty,\beta}[\mathbf{a},\mu]$   is based on the same argument). Hence, by the equalities in \eqref{jumpeqn2} and \eqref{jumpeqn3} one deduces the validity of \eqref{jumpeqn1}.\qed\medskip

\addcontentsline{toc}{chapter}{\bibname}
\bibliographystyle{plain}

\def\cprime{$'$} \def\cprime{$'$} \def\cprime{$'$}

\end{document}